\documentclass{article}

\usepackage{PRIMEarxiv}

\usepackage[utf8]{inputenc} % allow utf-8 input
\usepackage[T1]{fontenc}    % use 8-bit T1 fonts
\usepackage{hyperref}       % hyperlinks
\usepackage{url}            % simple URL typesetting
\usepackage{booktabs}       % professional-quality tables
\usepackage{amsfonts}       % blackboard math symbols
\usepackage{nicefrac}       % compact symbols for 1/2, etc.
\usepackage{microtype}      % microtypography
\usepackage{lipsum}
\usepackage{fancyhdr}       % header
\usepackage{graphicx}       % graphics
\usepackage{amsmath}
\usepackage{subcaption}
\graphicspath{{media/}}     % organize your images and other figures under media/ folder

\title{Stability analysis through folds: An end-loaded elastica with a lever arm}
\author{
  Siva Prasad Chakri Dhanakoti \\
  Department of Mathematics and Computer Science \\
  Freie Universität Berlin, Berlin 14135, Germany\\
  \texttt{chakri.dhanakoti@gmail.com} \\
  \AND
   Department of Civil, Environmental and Mechanical engineering \\
  University of Trento, Trento 38122, Italy\\
  \texttt{siva.dhanakoti@unitn.it} \\
}

\begin{document}
\maketitle

\begin{abstract}

Many physical systems can be modelled as parameter-dependent variational problems. In numerous cases, multiple equilibria co-exist, requiring the evaluation of their stability, and the monitoring of transitions between them. Generally, the stability characteristics of the equilibria change near folds in the parameter space. The direction of stability changes is embedded in a specific projection of the solutions, known as distinguished bifurcation diagrams. In this article, we identify such projections for variational problems characterized by fixed-free ends- a class of problems frequently encountered in mechanics. Using these diagrams, we study an Elastica subject to an end load applied through a rigid lever arm. Several instances of snap-back instability are reported, along with their dependence on system parameters through numerical examples. These findings have potential applications in the design of soft robot arms and other actuator designs. 

%( crucial for the physical understanding of the mechanical behavior of slender elastic rods. The equilibrium equations on their own are not sufficient to determine the local stability of the solutions; we also need to compute wether the solution is a local minimum or maximum of the system's energy.)

% Typically, one branch of solutions  in a bifurcation diagram represents stable solutions and other branch represents unstable solutions. We employ this information to examine the stability of variational systems which are characterized by fixed-free ends.

%Para 1 The abstract text goes here.The abstract text goes here. The abstract text goes here. The abstract text goes here.
%The abstract text goes here. The abstract text goes here. The abstract text goes here.\absbreak The abstract text goes here.
\end{abstract}

% keywords can be removed

\section{Introduction}
The study of physical problems using the variational framework has a rich history. In this setting, the equilibria are characterized as critical points of an energy-like functional. Typically, these problems are expressed in terms of a control parameter, generating a parameter family of equilibria. Indeed, the equilibria are computed by a slow scan of a parameter from a known solution using a continuation technique~\cite{DOEDEL1991a}. However, not all equilibria are stable and physically relevant; they must correspond to the local minima of the functional. Hence, it is crucial to determine whether the critical points are local minima (or stable) and track their changes. Generally, the stability of parameter-dependent equilibria changes at singularities, such as a fold or bifurcation in the parameter space~\cite{Sattinger1972, Rabinowitz1973, iooss2012elementary}. However, it remains unclear whether a fold in the unstable equilibrium results in a stable equilibrium or a higher-order unstable equilibrium. The direction of stability transitions can be encoded through a particular projection of solutions, using distinguished bifurcation diagrams, and the stability can be assessed with a little observation and without any additional analysis. This idea of topologically exploiting the variational structure for stability examination was first employed in Astrophysics~\cite{Katz1,Katz2}. Distinguished bifurcation diagrams were later developed for finite-dimensional mechanical systems~\cite{thompson1979stability} that are subsequently extended to infinite-dimensional systems~\cite{Maddocks1987, Hoffman2005}. Since then, several versions of these diagrams have appeared in investigations of elastic problems~\cite{neukirch2002a, Heijden2003}. In most studies, the diagrams are limited to cases with fixed-fixed ends or cases where the varying parameter appears within the integrand or as a Lagrange multiplier. In this study, we generalize the distinguished bifurcation diagrams to variational problems subject to fixed-free ends, when the bifurcation parameter appears in the boundary conditions. Consequently, two relevant cases arise: when the parameter appears at the fixed end and when it appears at the free end. For this analysis, we introduce the notion of a stability index that is equal to the dimension of the subspace over which the associated second variation operator is negative. The critical points with index zero are local minima and stable. Alternatively, the energy of the system can be used to determine the direction of stability~\cite{Maddocks1995,Oreilly2011}. The lower branch of the fold corresponds to the lower-energy state and therefore, local minima. 

The developed framework is applied to test the stability of elastica under an end load through a lever arm. Cantilever structures with a lever arm are common in nature and engineering—for example, in plants supporting fruits or in load-bearing flexible systems. In recent years, nonlinear rod models~\cite{antman2006nonlinear} have been actively employed in the fields of biophysics~\cite{Manning_DNA}, computer graphics~\cite{romero2021physical, Oreilly2012} and soft robot manipulators~\cite{Armanini2017, nayak2019}. Elastica, a simple planar version of this rod model, has equally received attention. A flexible elastic structure subjected to loads or constraints may exhibit multiple equilibria. Hence, it is crucial to determine their stability, as only stable structures are physically viable. Unstable equilibria may exist under specific constraints or boundary conditions and transitions to the available stable equilibrium under disturbances. Indeed, the presence of instability also hints at snap-back phenomena, a characteristic of multi-stable systems~\cite{Hu_2015}. Numerous studies have been conducted on the stability of elastic rods, which involve either the Jacobi test~\cite{Manning1998, Hoffman2002, Hoffman2004, Batista2015, dhanakoti2024stabilitycantileverlikestructuresapplications, Levyakov2010} or Hessian eigenvalue determination~\cite{goriely1997nonlinear, kumar2010generalized}. Alternative techniques based on bifurcation analysis are also available~\cite{Wang2021}. In this article, we assess the stability without performing any such rigorous computations. But, instead, conclusions are drawn by qualitatively examining the plot of appropriate projection of the equilibria against the parameter.

Advances in material science produced highly deformable alloys and polymers, accelerating the development of soft robotics. Generally, soft robots employ very flexible structures to generate compliant mechanisms. A flexible elastic structure under load can exhibit multiple equilibria, where one equilibrium may transition to another through snap-back instability under disturbances. These bistable or multi-stable systems exhibit large rapid displacements for small stimuli when operated around the regions of instabilities. As the system surpasses the initial energy barrier, it transitions to a new state, accompanied by a release of energy. Some of nature's quickest mechanisms seen in Mantis shrimp~\cite{patek2004deadly}, Venus fly trap~\cite{skotheim2005physical}, and Hummingbird beak~\cite{smith2011} can be attributed to this snap-back instability. These mechanisms are widely used in the design of high-performance soft actuators~\cite{luo2024intrinsically}, soft robots~\cite{Jin2023}, deployable and morphing structures~\cite{Yang2023} and haptic devices~\cite{Chen2018}. Inspired by these expanding applications, we aim to investigate elastica with a lever arm and the dependence of its stability on system parameters. These findings could offer invaluable insights for designing innovative mechanisms. Particularly, this study signifies the distinguished bifurcation diagrams for the analysis of multi-stable systems in engineering applications.

This article is organized as follows. The section~\ref{sec:Sec1} introduces the unconstrained calculus of variation problem with a parameter in the boundary conditions and describes the associated distinguished bifurcation diagrams. Section~\ref{sec:Sec2} formulates the elastica with an end load acting through a lever arm. A bifurcation analysis of resulting equilibria is also performed. In section~\ref{sec:Sec3}, we present several examples investigating the impact of varying parameters on stability using distinguished bifurcation diagrams. Finally, section~\ref{sec:Summary} provides a summary and discussion of the results.

\section{The Unconstrained Variational Problem}
\label{sec:Sec1}
In this section, we extend the theory of distinguished bifurcation diagrams developed by~\cite{Maddocks1987} and~\cite{Hoffman2005} to variational problems subject to fixed-free ends. We consider both cases, where the bifurcation parameter appears at the fixed end and when it appears at the free end. In section~\ref{ssec:Problem}, the problem and the relevant notation are introduced. The distinguished diagram at a simple fold when the parameter appears at the fixed end is developed in section~\ref{ssec:Fixed_end} and when the parameter appears at the free end is developed in section~\ref{ssec:Free_end}.

\subsection{Problem Formulation}
\label{ssec:Problem}
Let $\boldsymbol \zeta: s \rightarrow \mathbb{R}^{p}$ be a vector-valued continuous  differentiable function for $s\in[0,l]$. Consider an unconstrained calculus of variations problem of the form
 \begin{align}
 \label{eqn:parameter_functional_2}
     J ( \boldsymbol \zeta , \xi_{0},\xi_{l}) = \int_{0}^{l} \mathcal{L}(\boldsymbol \zeta, \boldsymbol \zeta^{\prime},s)ds + B(\boldsymbol \zeta,\xi_{l},l), 
 \end{align}
subject to a parameter-dependent \emph{Dirichlet  boundary condition} (or fixed boundary condition) at the end $s=0$
\begin{equation}
\boldsymbol \zeta(0)=\boldsymbol \zeta_{o}(\xi_{0}),
\end{equation}
while the other end $s=l$ is held free. The scalar function $\mathcal{L}$ has continuous second derivatives with respect to all its arguments and is convex in its second argument. On the other hand, the scalar function $B$ has continuous second derivatives with respect to its first argument and has continuous first derivatives with respect to the other arguments. The system can be influenced by the parameters $\xi_{0}$ and $\xi_{l}$ at the boundaries $s=0$ and $s=l$ respectively. In the current study, we vary either $\xi_{0}$ or $\xi_{l}$, but not both at a time. Maddocks~\cite{Maddocks1987} generated the bifurcation diagrams for problems when the parameter appears in the integrand and with homogeneous boundary conditions, which were later extended to the case when the parameter appears in the fixed end in a problem with  non-homogeneous boundary conditions~\cite{Hoffman2005}. However, these studies are confined to cases of fixed-fixed ends. Here, we extend the analysis to the case where one of the fixed ends is set free. The first-order condition yields the critical points $\boldsymbol \zeta_{o}(s,\tau,\xi_{i})$, $i=0$ or $l$, of~\eqref{eqn:parameter_functional_2} as solutions to the well-known \emph{Euler-Lagrange} equations 
\begin{align}
    \label{eqn:Euler_Lagrange_bif_dist}
    - \left( \frac{\partial \mathcal{L}}{\partial \boldsymbol \zeta^{\prime}}\right)^{\prime}
    + \frac{\partial \mathcal{L}}{\partial \boldsymbol \zeta} =\mathbf{0}, \qquad
    \boldsymbol \zeta(0)=\boldsymbol \zeta_{o}(\xi_{0}),
    \end{align}
    along with the natural boundary condition at the free end $s=l$
    \begin{align}
    \label{eqn:Natural_BCs}
    \left(\frac{\partial \mathcal{L}}{\partial \boldsymbol \zeta^{\prime}} + \frac{\partial B}{\partial \boldsymbol \zeta} \right)_{s=l}=\mathbf{0}.
 \end{align}
The notation $\xi_{i}$ refers to either $\xi_{0}$ or $\xi_{l}$ when not explicitly specified. These critical points $\boldsymbol \zeta_{o}(s,\tau,\xi_{i})$ exist as a continuous curve of solutions along a branch for some pseudo-arclength parameterization $\tau$~\cite{DOEDEL1991a} satisfying the property
\begin{align*}
 \big| \dot{ \boldsymbol \zeta}(s,\tau) \big| ^{2} + \big| \dot{ \xi_{i}}(\tau) \big| ^{2}=1.
\end{align*}
The notation $\dot{\boldsymbol \zeta}( s,\tau)$ denotes the derivative with respect to pseudo-arclength $\tau$, whereas ${\boldsymbol \zeta}^{\prime}( s,\tau)$ denotes the derivative with respect to independent parameter $s$. The second variation $\delta^{2}J$ evaluated at the critical points $\boldsymbol \zeta_{o}(s,\tau,\xi_{i})$, must satisfy the condition~\cite{gelfand1963calculus}:
 \begin{align}
\label{eqn:Second_order_cond} 
\begin{split}
   \delta ^{2}J(\boldsymbol \zeta_{o})[\mathbf{h}] =\frac{1}{2} & \int_{0}^{l}\left( \mathbf{h}^{\prime} \cdot \mathbf{P} \mathbf{h}^{\prime} +  \mathbf{h} \cdot \mathbf{C} \mathbf{h}^{\prime}+  \mathbf{h}^{\prime}  \cdot \mathbf{C}^{T} \mathbf{h}+ \mathbf{h} \cdot  \mathbf{Q} \mathbf{h} \right)ds + \frac{1}{2}\mathbf{B} \mathbf{h}(l)   \cdot   \mathbf{h}(l)  \geq 0,
   \end{split}
  \end{align}
so that they correspond to the minima of the functional $J$, where $\mathbf{P},\mathbf{C}$, $\mathbf{Q}$ and $\mathbf{B}$ are $p \times p$ Hessian matrices evaluated at the extremal $\boldsymbol \zeta_{o}$ given by
  \begin{align}
  \label{eqn:PQC_matrices}
    \begin{split}
      \mathbf{P}=\mathcal{L}_{ \boldsymbol \zeta^{\prime}  \boldsymbol \zeta^{\prime}}( \boldsymbol \zeta_{o},  \boldsymbol \zeta_{o}^{\prime},s),  \quad \mathbf{C} =\mathcal{L}_{ \boldsymbol \zeta  \boldsymbol\zeta^{\prime} }( \boldsymbol \zeta_{o},  \boldsymbol \zeta_{o}^{\prime},s),\quad 
      \mathbf{Q}=\mathcal{L}_{ \boldsymbol \zeta  \boldsymbol \zeta}( \boldsymbol \zeta_{o},  \boldsymbol \zeta_{o}^{\prime},s), \quad \text{and} \quad \mathbf{B}=B_{\boldsymbol \zeta \boldsymbol \zeta}(\boldsymbol \zeta_{o},l).
      \end{split}
  \end{align}
  For brevity, the dependence of $\mathbf{P},\mathbf{C}$ and $\mathbf{Q}$ on $s$ is not explicitly shown. The matrices $\mathbf{P}$,$\mathbf{Q}$ and $\mathbf{B}$ are symmetric, whereas the matrix $\mathbf{C}$ may not be. Here, $\mathbf{h}$ is a variation in the solution that satisfy the linearized boundary condition:
  \begin{align}
  \mathbf{h}(0)=\mathbf{0}.
  \end{align}
The integrand in~\eqref{eqn:Second_order_cond} has no explicit dependence on $\tau$ and $\xi_{i}$, and the coefficient matrices depend on them through the solutions $\boldsymbol \zeta_{o}(s, \xi_{i})$ and $\boldsymbol \zeta_{o}^{\prime}(s, \xi_{i})$. We also assume that the Legendre's strengthened condition holds
    \begin{align}
            \mathbf{P} > 0,
    \end{align}
    i.e., the symmetric matrix $\mathbf{P}$ is positive definite.
    The application of integration by parts on~\eqref{eqn:Second_order_cond} and the vanishing boundary term leads to an alternate form of the second variation
    \begin{align}
    \label{eqn:secondvar}
    \delta^{2}J[\mathbf{h}] =\frac{1}{2}\langle \mathcal{S}\mathbf{h} , \mathbf{h}\rangle +  \frac{1}{2}\left( \mathbf{P} \mathbf{h}^{\prime} +\mathbf{C}^{T}\mathbf{h} + \mathbf{B}\mathbf{h}\right)\cdot \mathbf{h} \bigg|_{s=l} , 
  \end{align}
  where $\mathcal{S}$ is the second-order differential operator
  \begin{align}
  \label{eqn:accessory_bvp}
      \mathcal{S} \mathbf{h} \equiv - \frac{d}{ds} \left( \mathbf{P} \mathbf{h}^{\prime} +\mathbf{C}^{T}\mathbf{h} \right) + \mathbf{C}  \mathbf{h}^{\prime} + \mathbf{Q} \mathbf{h} , 
  \end{align}
   and $\langle \cdot , \cdot \rangle$ represents the standard inner-product in $L^{2}$- space. The coefficient $\frac{1}{2}$ has no effect on the results and is omitted in the subsequent analysis. We define an notion of index~\cite{Morse1951} which denotes the dimension of the $\mathbf{h}$ subspace over which the second variation denoted by~\eqref{eqn:secondvar} is negative. The zero index solutions satisfy the necessary condition for local minima  given by~\eqref{eqn:Second_order_cond} and correspond to stable solutions. In fact, the index can be interpreted  as the number of eigenvalues of the eigenvalue problem
   \begin{align}
   \label{eqn:eigen_valueproblem}
\mathcal{S} \mathbf{h}= \mu  \mathbf{h}, \qquad \mathbf{h}(0)=\mathbf{0},
\end{align}
that generates negative
\begin{align}
\label{eqn:second_variation}
\delta^{2}J[\mathbf{h}] =\mu\langle \mathbf{h}, \mathbf{h}\rangle  +  \left( \mathbf{P} \mathbf{h}^{\prime} +\mathbf{C}^{T}\mathbf{h} + \mathbf{B}\mathbf{h}\right)\cdot \mathbf{h} \bigg|_{s=l}. 
\end{align}

The operator $\mathcal{S}$ is assumed to have eigenvalues and eigenfunctions that smoothly depend on the pseudo-arclength $\tau$. Therefore, it is not sufficient to merely determine the sign of the eigenvalue; it is also necessary to account for the additional boundary contribution. Alternatively, one can consider the~\eqref{eqn:Second_order_cond} directly using numerical methods. Upon discretization, it leads to a Global Stiffness Matrix whose eigenvalues dictate stability - positive eigenvalues imply local minimum and thus indicate stability. The index for the case of fixed-fixed ends is shown to be equal to the number of conjugate points (for more details refer to~\cite{Manning1998}). However, we will not pursue our analysis in that direction. Instead, we track changes in the index at critical points. 

We now state an integration by parts result for the $\mathcal{S}$ operator, which is used frequently in the subsequent analysis. For any two vector-valued continuous functions $\mathbf{u},\mathbf{v}$, we can easily prove 
\begin{align}
\label{eqn:lma1}
 \langle \mathcal{S}\mathbf{u},\mathbf{v}\rangle= \langle \mathbf{u},\mathcal{S}\mathbf{v}\rangle + \Big[ \left( \mathbf{P} \mathbf{v}^{\prime} +\mathbf{C}^{T}\mathbf{v} \right) \cdot \mathbf{u} \Big]^{l}_{0} -  \Big[ \left( \mathbf{P} \mathbf{u}^{\prime} +\mathbf{C}^{T}\mathbf{u} \right) \cdot \mathbf{v} \Big]^{l}_{0},
\end{align} 
 through successive applications of integration by parts. Therefore, when both the boundary terms vanish $\mathcal{S}$ is a self-adjoint operator.

 Parameter-dependent variational problems generate a family of extremals, and in the case of a single parameter $\xi$, a curve of solutions is obtained. Sometimes, these solutions are characterized by \emph{folds}- points at which the curve slopes vertically. In bifurcation terminology, these points are also referred to as \emph{saddle-node bifurcation} or \emph{turning points}. On differentiating the Euler-Lagrange equations given by~\eqref{eqn:Euler_Lagrange_bif_dist} and the non-linear boundary conditions given by~\eqref{eqn:Natural_BCs} with respect to pseudo-arclength $\tau$, we obtain the following expressions for the two cases of varying parameters
\begin{subequations}
   \label{eqn:Derivative_EL_equations}
  \begin{align}
      \mathcal{S}\dot{\boldsymbol \zeta} &=\mathbf{0},  
     \qquad  \dot{\boldsymbol \zeta}(0)= \frac{d \boldsymbol  \zeta_{o}}{d \xi} \dot{\xi}_{0},  
   \quad \mathbf{P} \dot{\boldsymbol \zeta^{\prime}}(l)  + \mathbf{C}^{T}\dot{\boldsymbol \zeta}(l)  + \mathbf{B}\dot{\boldsymbol \zeta}(l) =\mathbf{0}, \\
     \mathcal{S}\dot{\boldsymbol \zeta} &=\mathbf{0},   \qquad \dot{\boldsymbol \zeta}(0)=\mathbf{0},
   \quad \qquad  \mathbf{P} \dot{\boldsymbol \zeta^{\prime}}(l) + \mathbf{C}^{T}\dot{\boldsymbol \zeta}(l) + \mathbf{B} \dot{\boldsymbol \zeta}(l)+ \frac{\partial^{2} B}{\partial \boldsymbol \zeta \partial \xi_{l}} \dot{\xi}_{l} =\mathbf{0}.
  \end{align}
\end{subequations}
The Hessian matrices that appear in the boundary terms $s=l$ are evaluated at that boundary, and their dependence on $l$ is not included for conciseness. Each of these boundary value problems corresponds to the eigenvalue problem denoted by~\eqref{eqn:eigen_valueproblem} with zero eigenvalue whenever $\dot{\xi}_{0}=0$ or $\dot{\xi}_{l}=0$, with the eigenvector $\mathbf{h}$ coinciding with $\dot{\boldsymbol \zeta}$. The boundary condition at $s=l$ also removes the boundary term in~\eqref{eqn:second_variation}, fetching a zero second-variation quadratic functional. We focus on the behavior of the motion of this functional around these special points - folds. A \emph{simple fold} further satisfies:
\begin{subequations}
 \begin{align}
&\ddot{\xi}_{i} \neq 0, \qquad \\
\label{eqn:transversal}
&\frac{d}{d \tau} J_{\xi_{i}} \neq 0, \qquad \\
 &\text{and } \mu=0\text{ is a simple eigenvalue of }\mathcal{S}.
\end{align}
\end{subequations}
An eigenvalue is called a simple eigenvalue if its algebraic multiplicity is one. The meaning of~\eqref{eqn:transversal} becomes clearer towards the end. We analyze the behavior of the second-variational functional in the vicinity of these folds, which involves simultaneously monitoring the eigenvalues and the boundary term contribution from the corresponding eigenfunctions. We assume that changes in stability occur exclusively in the vicinity of critical points, such as folds or bifurcation points~\cite{Sattinger1972}.

\iffalse
Let $\mathbf{h}$ denote the eigenvector corresponding to the eigenvalue $\mu$.

On taking the inner product of~\eqref{eqn:Derivative_EL_equations} with the eigenvector $\mathbf{h}$, we obtain
\begin{align}
\langle  \mathcal{S} \dot{\boldsymbol \zeta}, \mathbf{h} \rangle = 0.
\end{align}
Integrating by parts fetches 
\begin{align}
\label{eqn:IBP_Derivative}
    \langle   \mathcal{S} \dot{\boldsymbol \zeta} , \mathbf{h} \rangle = \langle  \mathcal{S} \mathbf{h} , \dot{\boldsymbol \zeta} \rangle  - \left[ \left(\mathbf{P}  \mathbf{h}^{\prime} + \mathbf{C}^{T}  \mathbf{h}  + \mathbf{B}  \mathbf{h}\right) \cdot \dot{\boldsymbol \zeta} \right]_{0}^{l} + \left[ \left(\mathbf{P} \dot{\boldsymbol \zeta}^{\prime} + \mathbf{C}^{T} \dot{\boldsymbol \zeta} +  \mathbf{B} \dot{\boldsymbol \zeta }\right) \cdot  \mathbf{h} \right]_{0}^{l}=0.
\end{align}
 The key aspect of the analysis is that, in the vicinity of a fold, an elementary expression can be found for the derivative of the critical eigenvalue $\dot{\mu}$.

\begin{align}
\langle \mathcal{S}\mathbf{h}, \mathbf{h}\rangle +\left( \mathbf{P} \mathbf{h}^{\prime} +\mathbf{C}^{T}\mathbf{h} + \mathbf{B}\mathbf{h}\right)\cdot \mathbf{h} \bigg|_{s=l} 
\end{align}
We are interested in tracking its direction in the neighborhood of a fold, where it has a zero value.
\fi

On differentiating the second variational quadratic form denoted by~\eqref{eqn:secondvar} with respect to $\tau$, we obtain:
\begin{align}
\label{eqn:second_var_dot}
\frac{d}{d \tau} \delta^{2}J[\mathbf{h}]=\langle \dot{ \left(\mathcal{S}\mathbf{h} \right) }, \mathbf{h}\rangle+ \dot{\bigg( \mathbf{P} \mathbf{h}^{\prime} +\mathbf{C}^{T}\mathbf{h} + \mathbf{B}\mathbf{h}\bigg)}\cdot \mathbf{h} \bigg|_{s=l} + \langle  \mathcal{S}\mathbf{h}  , \dot{\mathbf{h}} \rangle +\bigg( \mathbf{P} \mathbf{h}^{\prime} +\mathbf{C}^{T}\mathbf{h} + \mathbf{B}\mathbf{h}\bigg)\cdot \dot{\mathbf{h}} \bigg|_{s=l}, 
\end{align}
along with the boundary condition
\begin{align}
\label{eqn:h_dot}
\dot{\mathbf{h}}(0)=\mathbf{0}.
\end{align}
Near the fold where $\mathbf{h}\equiv \dot{\boldsymbol \zeta}$ ,~\eqref{eqn:second_var_dot} reduces to:
\begin{align*}
\frac{d}{d \tau} \delta^{2}J[\mathbf{h}]=\langle \frac{d}{d \tau} \mathcal{S}\mathbf{h}  , \mathbf{h}\rangle+ \mathbf{h}(l) \cdot \frac{d}{d \tau}{\left( \mathbf{P} \mathbf{h}^{\prime}(l) +\mathbf{C}^{T}\mathbf{h}(l) + \mathbf{B}\mathbf{h}(l)\right)}.
\end{align*}
Upon further expansion and replacing $\mathbf{h}$ with $\dot{\boldsymbol \zeta}$, we get
\begin{align}
\label{eqn:del2J_dot}
\frac{d}{d \tau} \delta^{2}J[\mathbf{h}]= \langle  \mathcal{S}\dot{\mathbf{h}} , \dot{\boldsymbol \zeta}\rangle+
\langle \dot{ \mathcal{S}} \dot{\boldsymbol \zeta} , \dot{\boldsymbol \zeta}\rangle+  \left( \dot{\mathbf{P} }\dot{\boldsymbol \zeta}^{\prime} + \dot{\mathbf{C}}^{T}\dot{\boldsymbol \zeta} + \dot{\mathbf{B}}\dot{\boldsymbol \zeta}\right) \cdot \dot{\boldsymbol \zeta}  \bigg|_{s=l} + \left( \mathbf{P} \dot{\mathbf{h}}^{\prime} + \mathbf{C}^{T}\dot{\mathbf{h}} + \mathbf{B} \dot{\mathbf{h}}\right) \cdot \dot{\boldsymbol \zeta}  \bigg|_{s=l}.
\end{align}
Here, the operator $\dot{\mathcal{S}}$ denotes the derivative of $\mathcal{S}$ with respect to $\tau$ that takes the form
\begin{align*}
\dot{\mathcal{S}} \mathbf{h} = \left(\frac{d}{d \tau} \mathcal{S}\right) \mathbf{h}=  - \frac{d}{ds} \left( \dot{\mathbf{P}} \mathbf{h}^{\prime} +\dot{\mathbf{C}}^{T}\mathbf{h} \right) + \dot{\mathbf{C}  }\mathbf{h}^{\prime} + \dot{\mathbf{Q}} \mathbf{h},
\end{align*}
and $\dot{\mathbf{P}},\dot{\mathbf{C}},\dot{\mathbf{Q}}$ and $\dot{\mathbf{B}}$ are the derivatives of the Hessian matrices with respect to $\tau$. Consider the term $\langle  \mathcal{S}\dot{\mathbf{h}} , \dot{\boldsymbol \zeta}\rangle$. Applying the integration by parts identity given by~\eqref{eqn:lma1}, we obtain
\begin{align*}
\langle  \mathcal{S}\dot{\mathbf{h}} , \dot{\boldsymbol \zeta}\rangle= \langle  \dot{\mathbf{h}} , \mathcal{S}\dot{\boldsymbol \zeta}\rangle - \left[\left( \mathbf{P} \dot{\mathbf{h}}^{\prime} + \mathbf{C}^{T}\dot{\mathbf{h}} + \mathbf{B}\dot{\mathbf{h}}\right) \cdot \dot{\boldsymbol \zeta} \right]_{0}^{l} + \left[\left( \mathbf{P} \dot{\boldsymbol \zeta}^{\prime} + \mathbf{C}^{T}\dot{\boldsymbol \zeta} + \mathbf{B}\dot{\boldsymbol \zeta}\right) \cdot \dot{\mathbf{h}} \right]_{0}^{l}
\end{align*}
Since $\dot{\boldsymbol \zeta}$ is the zero eigenvector of $\mathcal{S}$, we have $\mathcal{S} \dot{\boldsymbol \zeta}=\mathbf{0}$ implying $\langle  \dot{\mathbf{h}} , \mathcal{S}\dot{\boldsymbol \zeta}\rangle=0$. On applying the boundary conditions denoted by~\eqref{eqn:Derivative_EL_equations} and~\eqref{eqn:h_dot} at a fold $\dot{\xi_{i}}$, the above expression simplifies to
\begin{align}
\label{eqn:Sh_dot}
\langle  \mathcal{S}\dot{\mathbf{h}} , \dot{\boldsymbol \zeta}\rangle= -\left(\mathbf{P} \dot{\mathbf{h}}^{\prime} + \mathbf{C}^{T}\dot{\mathbf{h}} + \mathbf{B}\dot{\mathbf{h}}\right) \cdot \dot{\boldsymbol \zeta} \bigg|_{s=l}.
\end{align}
 Now consider the second term $ \langle \dot{\mathcal{S}} \dot{\boldsymbol \zeta},\dot{\boldsymbol \zeta} \rangle$. Differentiating~\eqref{eqn:Derivative_EL_equations} with respect to $\tau$ gives the identity $ \dot{\mathcal{S}} \dot{\boldsymbol \zeta} = -\mathcal{S} \ddot{\boldsymbol \zeta}$
along with the following sets of boundary conditions. When the varying parameter appears at the fixed end $s=0$, we have
\begin{subequations}
\begin{align}
\ddot{\boldsymbol \zeta}(0)&= \frac{\partial^{2} \boldsymbol \zeta_{o}}{ \partial \xi_{0}^{2}} \dot{\xi}_{0}^{2} + \frac{\partial \boldsymbol \zeta_{o}}{\partial \xi_{0}} \ddot{\xi}_{0} =  \frac{\partial \boldsymbol \zeta_{o}}{ \partial \xi_{0}} \ddot{\xi}_{0}, \\
 \mathbf{P} \ddot{\boldsymbol \zeta}(l) &+\mathbf{C}^{T}  \ddot{\boldsymbol \zeta^{\prime}}(l) +\mathbf{B}  \ddot{\boldsymbol \zeta^{\prime}}(l)  + \dot{\mathbf{P}} \dot{\boldsymbol \zeta}(l) + \dot{\mathbf{C}^{T} } \dot{\boldsymbol \zeta^{\prime}}(l) + \dot{\mathbf{B} } \dot{\boldsymbol \zeta^{\prime}}(l)=\mathbf{0},
\label{eqn:boundaryterms_fixedend_b}
\end{align}
\label{eqn:boundaryterms_fixedend}
\end{subequations}
and when the varying parameter appears at the fixed end $s=l$, we have
\begin{subequations}
\begin{align}
\ddot{\boldsymbol \zeta}(0)&= \mathbf{0}, \\
\begin{split}
 \mathbf{P} \ddot{\boldsymbol \zeta} (l)&+\mathbf{C}^{T}  \ddot{\boldsymbol \zeta^{\prime}}(l) +\mathbf{B}  \ddot{\boldsymbol \zeta^{\prime}}(l)  + \dot{\mathbf{P}} \dot{\boldsymbol \zeta}(l) + \dot{\mathbf{C}}^{T}  \dot{\boldsymbol \zeta^{\prime}}(l) + \dot{\mathbf{B}} \dot{\boldsymbol \zeta^{\prime}}(l)\\& +\frac{\partial^{3} B}{\partial \boldsymbol \zeta \partial \xi_{l}^{2}} \dot{\xi}_{l}^{2}+ \frac{\partial^{3} B}{\partial \boldsymbol \zeta^{2} \partial \xi_{l}} \dot{\boldsymbol{\zeta}}(l) \dot{\xi}_{l}^{2} +\frac{\partial^{2} B}{\partial \boldsymbol \zeta \partial \xi_{l}} \ddot{\xi}_{l}=\mathbf{0}, 
\end{split}
\label{eqn:boundaryterms_freeend_b}
\end{align}
\label{eqn:boundaryterms_freeend}
\end{subequations}
Consequently, the term $- \langle \mathcal{S} \ddot{\boldsymbol \zeta}, \dot{\boldsymbol \zeta}\rangle$ can be rewritten using the identity given by~\eqref{eqn:lma1}) as  
\begin{align*}
     - \langle \mathcal{S} \ddot{\boldsymbol \zeta}, \dot{\boldsymbol \zeta} \rangle& = - \langle \ddot{\boldsymbol \zeta} ,\mathcal{S}  \dot{\boldsymbol \zeta}\rangle  + \bigg[\left(\mathbf{P} \ddot{\boldsymbol \zeta}^{\prime} + \mathbf{C}^{T}\ddot{\boldsymbol \zeta} + \mathbf{B}\ddot{\boldsymbol \zeta} \right) \cdot \dot{\boldsymbol \zeta} \bigg]_{0}^{l} - \bigg[ \left(\mathbf{P} \dot{\boldsymbol \zeta}^{\prime}+ \mathbf{C}^{T}\dot{\boldsymbol \zeta} + \mathbf{B}\dot{\boldsymbol \zeta} \right) \cdot \ddot{\boldsymbol \zeta } \bigg]_{0}^{l}.
 \end{align*}
Since $\mathcal{S} \dot{\boldsymbol \zeta}=\mathbf{0}$, this expression simplifies to
 \begin{align}
 \label{eqn:S_ddot}
     - \langle \mathcal{S} \ddot{\boldsymbol \zeta}, \dot{\boldsymbol \zeta} \rangle& = \bigg[\left(\mathbf{P} \ddot{\boldsymbol \zeta}^{\prime} + \mathbf{C}^{T}\ddot{\boldsymbol \zeta} + \mathbf{B}\ddot{\boldsymbol \zeta} \right) \cdot \dot{\boldsymbol \zeta} \bigg]_{0}^{l} - \bigg[ \left(\mathbf{P} \dot{\boldsymbol \zeta}^{\prime}+ \mathbf{C}^{T}\dot{\boldsymbol \zeta} + \mathbf{B}\dot{\boldsymbol \zeta} \right) \cdot \ddot{\boldsymbol \zeta } \bigg]_{0}^{l}.
 \end{align}
From this point, we separate the analyses for the cases of varying parameter at the fixed end $(\xi_{0})$ and the free end $(\xi_{l})$.

\subsection{Varying Parameter at the Fixed-end}
\label{ssec:Fixed_end}
In this section, we analyze the behavior of $- \langle \mathcal{S} \ddot{\boldsymbol \zeta}, \dot{\boldsymbol \zeta} \rangle$ and consequently $\frac{d}{d \tau} \delta^{2}J$ when $\xi_{0}$ is varied. At the fold $\dot{\xi}_{0}=0$,~\eqref{eqn:Derivative_EL_equations} leads to $\dot{\boldsymbol \zeta}(0)=\mathbf{0}$ and $\left( \mathbf{P} \dot{\boldsymbol \zeta}^{\prime}+ \mathbf{C}^{T}\dot{\boldsymbol \zeta}+ \mathbf{B}\dot{\boldsymbol \zeta}\right)_{s=l}=\mathbf{0} $. Then,~\eqref{eqn:S_ddot} simplifies to  
 \begin{align*}
- \langle \mathcal{S} \ddot{\boldsymbol \zeta}, \dot{\boldsymbol \zeta} \rangle& =\left(\mathbf{P} \dot{\boldsymbol \zeta}^{\prime} + \mathbf{C}^{T}\dot{\boldsymbol \zeta} +  \mathbf{B}\dot{\boldsymbol \zeta} \right) \cdot \ddot{\boldsymbol \zeta} \bigg|_{s=0} +\left(\mathbf{P} \ddot{\boldsymbol \zeta}^{\prime} + \mathbf{C}^{T}\ddot{\boldsymbol \zeta} + \mathbf{B}\ddot{\boldsymbol \zeta} \right) \cdot \dot{\boldsymbol \zeta} \bigg|_{s=l}.
\end{align*}
Using the relations~\eqref{eqn:boundaryterms_fixedend} and recalling that $\mathbf{P}=\mathbf{P}^{T}$, we deduce 
\begin{align*}
- \langle \mathcal{S} \ddot{\boldsymbol \zeta}, \dot{\boldsymbol \zeta} \rangle&= \ddot{\xi}_{0} \left( \frac{\partial \boldsymbol \zeta_{o}}{\partial \xi_{0} }\cdot \mathbf{P} \dot{\boldsymbol \zeta}^{\prime}\right)\Bigg|_{s=0} + \left(\mathbf{C}^{T}\dot{\boldsymbol \zeta} +  \mathbf{B}\dot{\boldsymbol \zeta} \right) \cdot \ddot{\boldsymbol \zeta} \bigg|_{s=0} +\left(\mathbf{P} \ddot{\boldsymbol \zeta}^{\prime} + \mathbf{C}^{T}\ddot{\boldsymbol \zeta} + \mathbf{B}\ddot{\boldsymbol \zeta} \right) \cdot \dot{\boldsymbol \zeta} \bigg|_{s=l},\\
&=\ddot{\xi}_{0} \left( \frac{\partial \boldsymbol \zeta_{o}}{\partial \xi_{0} }\cdot \mathbf{P} \dot{\boldsymbol \zeta}^{\prime}\right)\Bigg|_{s=0}  +\left(\mathbf{P} \ddot{\boldsymbol \zeta}^{\prime} + \mathbf{C}^{T}\ddot{\boldsymbol \zeta} + \mathbf{B}\ddot{\boldsymbol \zeta} \right) \cdot \dot{\boldsymbol \zeta} \bigg|_{s=l}.
\end{align*}
The $\ddot{\boldsymbol \zeta}$ and $\ddot{\boldsymbol \zeta}^{\prime}$ terms can be eliminated using~\eqref{eqn:boundaryterms_fixedend_b}, resulting in 
\begin{align}
\label{eqn:s_ddot}
- \langle \mathcal{S} \ddot{\boldsymbol \zeta}, \dot{\boldsymbol \zeta} \rangle=\ddot{\xi}_{0} \left( \frac{\partial \boldsymbol \zeta_{o}}{\partial \xi_{0}} \cdot \mathbf{P} \dot{\boldsymbol \zeta}^{\prime}\right)\Bigg|_{s=0}  -\left(\dot{\mathbf{P}} \dot{\boldsymbol \zeta}^{\prime} + \dot{\mathbf{C}}^{T}\dot{\boldsymbol \zeta}  + \dot{\mathbf{B}} \dot{\boldsymbol \zeta} \right) \cdot \dot{\boldsymbol \zeta} \bigg|_{s=l}.
\end{align}
In the end, after substituting the~\eqref{eqn:Sh_dot} and~\eqref{eqn:s_ddot} in~\eqref{eqn:del2J_dot}, we obtain
 \begin{align}
 \label{eqn:deltasquare_dot}
 \frac{d}{d \tau}{\delta^{2}J}  =\ddot{\xi}_{0} \left( \frac{\partial \boldsymbol \zeta_{o}}{\partial \xi_{0}} \cdot\mathbf{P} \dot{\boldsymbol \zeta}^{\prime}\right)\Bigg|_{s=0}.
 \end{align}
 For the symmetric matrix $\mathbf{P} = \mathcal{L}_{\boldsymbol \zeta^{\prime}\boldsymbol \zeta^{\prime}}$, we can write
\begin{align*}
    \dot{\boldsymbol \zeta}^{\prime} \cdot \mathcal{L}_{\boldsymbol \zeta^{\prime} \boldsymbol \zeta^{\prime}} \frac{\partial \boldsymbol \zeta_{o}}{\partial \xi_{0}}\Bigg|_{s=0} &= \bigg(\dot{\boldsymbol \zeta}^{\prime} \cdot \mathcal{L}_{\boldsymbol \zeta^{\prime} \boldsymbol \zeta^{\prime }} \frac{\partial \boldsymbol \zeta_{o}}{\partial \xi_{0}} + \mathcal{L}_{\boldsymbol \zeta \boldsymbol \zeta^{\prime}} \frac{\partial \boldsymbol \zeta_{o}}{\partial \xi}\dot{\xi_{0}} \cdot \frac{\partial \boldsymbol \zeta_{o}}{\partial \xi_{0}}\bigg)\Bigg|_{s=0},  \\
    &= \bigg(\mathcal{L}_{\boldsymbol \zeta^{\prime} \boldsymbol \zeta^{\prime }}  \dot{\boldsymbol \zeta}^{\prime} \cdot \frac{\partial \boldsymbol \zeta_{o}}{\partial \xi_{0}} + \mathcal{L}_{\boldsymbol \zeta \boldsymbol \zeta^{\prime}} \frac{\partial \boldsymbol \zeta_{o}}{\partial \xi_{0}}\dot{\xi_{0}} \cdot \frac{\partial \boldsymbol \zeta_{o}}{\partial \xi_{0}}\bigg)\Bigg|_{s=0}=\frac{d }{d \tau} \bigg( \mathcal{L}_{\boldsymbol \zeta^{\prime}} \cdot \frac{\partial \boldsymbol \zeta_{o}}{\partial \xi_{0}}\bigg)\Bigg|_{s=0}.
\end{align*}
Then the expression for the derivative becomes 
\begin{align}
\label{eqn:distinguished_biff_ordinate_fixed}
    \frac{d}{d \tau} \delta^{2}J= \ddot{\xi_{0}} \frac{d}{d \tau} \bigg(\mathcal{L}_{\boldsymbol \zeta^{\prime}} \cdot \frac{\partial \boldsymbol \zeta_{o}}{\partial \xi_{0}}\bigg) \Bigg|_{s=0}.
\end{align}
For a simple fold, the terms $\ddot{\xi}_{0} $ and 
\begin{align}
\label{eqn:transversal_fixedend}
\frac{d}{d \tau} \bigg(\mathcal{L}_{\boldsymbol \zeta^{\prime}} \cdot \frac{\partial \boldsymbol \zeta_{o}}{\partial \xi_{0}}\bigg) \Bigg|_{s=0} 
\end{align}
are non-zero, and~\eqref{eqn:distinguished_biff_ordinate_fixed} provides the direction in which the second variation crosses zero. The non-zero value of~\eqref{eqn:transversal_fixedend} is same as the condition~\eqref{eqn:transversal}. The sign of $\ddot{\xi}_{0}$ indicates the direction in which the fold opens. A positive $\ddot{\xi}_{0}$ represents the fold opening towards the right, while a negative $\ddot{\xi}_{0}$ represents the fold opening toward the left. The distinguished bifurcation diagram is defined as the plot of $\left( \mathcal{L}_{\boldsymbol \zeta^{\prime}} \cdot \frac{\partial \boldsymbol \zeta_{o}}{\partial \xi_{0}}\right)\Bigg|_{s=0} $ vs. $\xi_{0}$ along the branches of critical points. As a simple fold is traversed, the critical eigenvalue $\mu$ and the associated eigenvector result in a sign change of $\delta^{2}J$ dictated by the signs of $\ddot{\xi}_{0}$ and $ \left( \mathcal{L}_{\boldsymbol \zeta^{\prime}} \cdot \frac{\partial \boldsymbol \zeta_{o}}{\partial \xi_{0}}\right)\Bigg|_{s=0}$. Therefore, the index, which is the number of directions along which the second variation is negative is altered by one at each fold, and the direction of change can be found by examining the shape of the projection of the solution branch as illustrated in Figure~\ref{fig:schematic_bifurcation_diagram_a}. Indeed, the expression for the ordinate is the same as in the case of fixed-fixed ends~\cite{Hoffman2005}, where the varying parameter appears at the fixed end $s=l$.

\subsection{Varying Parameter at the Free-end}
\label{ssec:Free_end}
\begin{figure}[t!]
    \centering
    \begin{subfigure}[t]{0.497\textwidth}
    \includegraphics[width=1.0\textwidth]{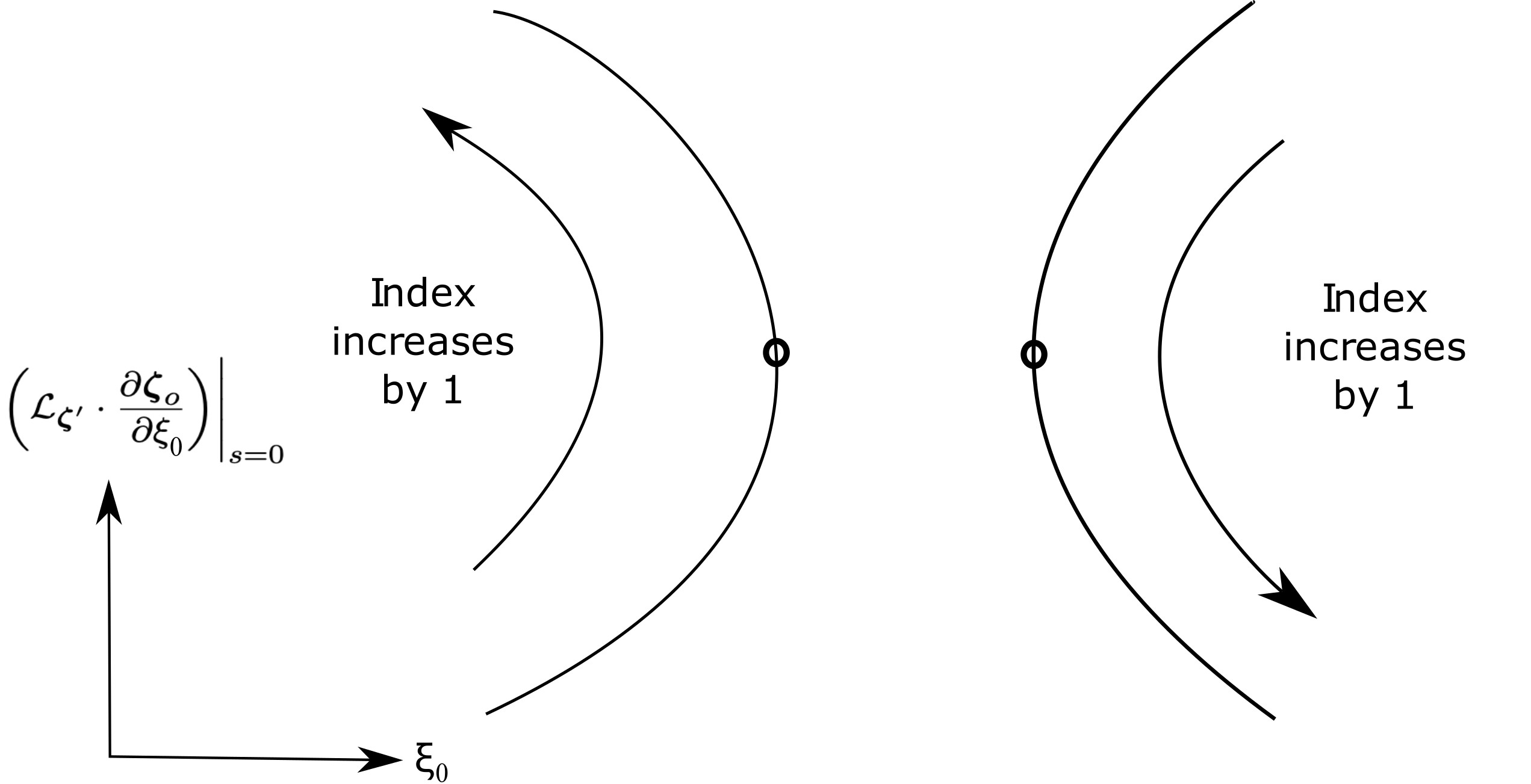}
    \caption{}
    \label{fig:schematic_bifurcation_diagram_a}
    \end{subfigure}
    \begin{subfigure}[t]{0.497\textwidth}
    \includegraphics[width=1.0\textwidth]{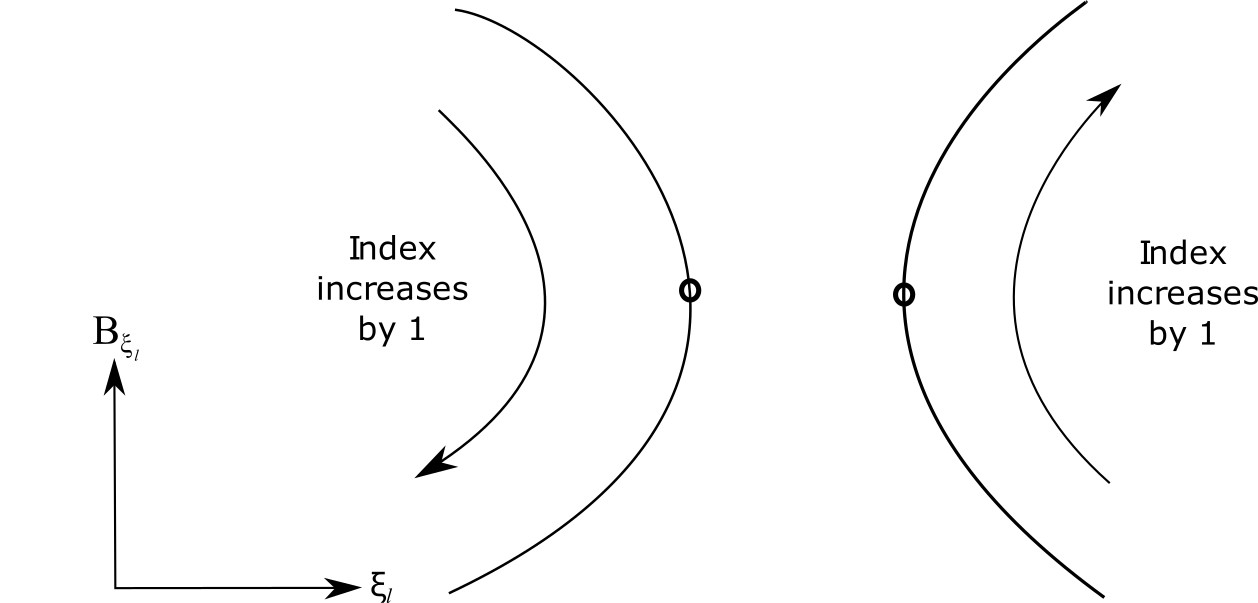}
    \caption{}
    \label{fig:schematic_bifurcation_diagram_b}
    \end{subfigure} 
    \caption{The direction of index change near simple folds in the distinguished bifurcation diagram for problems with fixed-free ends: a) Stability transition when the parameter in the fixed end $s=0$ is varied. b) Stability transition when the parameter in the free end $s=l$ is varied.}
    \label{fig:Distinguished_Bifurcation}
\end{figure}

Now consider the other case, where the parameter is at the free end $s=l$, namely $\xi_{l}$. At the fold $\dot{\xi}_{l}=0$, ~\eqref{eqn:Derivative_EL_equations} leads to $\dot{\boldsymbol \zeta}(0)=0$ and $\left( \mathbf{P} \dot{\boldsymbol \zeta}^{\prime}+ \mathbf{C}^{T}\dot{\boldsymbol \zeta}+ \mathbf{B}\dot{\boldsymbol \zeta}\right)_{s=l}=\mathbf{0}$. As a result, the expression for $ - \langle \mathcal{S} \ddot{\boldsymbol \zeta}, \dot{\boldsymbol \zeta} \rangle$ given by~\eqref{eqn:S_ddot} turns to
\begin{align*}
- \langle \mathcal{S} \ddot{\boldsymbol \zeta}, \dot{\boldsymbol \zeta} \rangle =\left(\mathbf{P} \ddot{\boldsymbol \zeta}^{\prime} + \mathbf{C}^{T}\ddot{\boldsymbol \zeta} + \mathbf{B}\ddot{\boldsymbol \zeta} \right) \cdot \dot{\boldsymbol \zeta} \bigg|_{l} +  \left(\mathbf{P} \dot{\boldsymbol \zeta}^{\prime}+ \mathbf{C}^{T}\dot{\boldsymbol \zeta} + \mathbf{B}\dot{\boldsymbol \zeta} \right) \cdot \ddot{\boldsymbol \zeta } \bigg|_{0}.
\end{align*}
Furthermore, the boundary conditions specified in~\eqref{eqn:boundaryterms_freeend_b} result in 
\begin{align}
\label{eqn:S_ddot2}
- \langle \mathcal{S} \ddot{\boldsymbol \zeta}, \dot{\boldsymbol \zeta} \rangle=-\left(\dot{\mathbf{P}} \ddot{\boldsymbol \zeta}^{\prime} + \dot{\mathbf{C}}^{T}\ddot{\boldsymbol \zeta} + \dot{\mathbf{B}}\ddot{\boldsymbol \zeta}+\frac{\partial^{2} B}{\partial \boldsymbol \zeta \partial \xi_{l}} \ddot{\xi}_{l} \right) \cdot \dot{\boldsymbol \zeta} \bigg|_{s=l}.
\end{align}
On substituting~\eqref{eqn:Sh_dot} and~\eqref{eqn:S_ddot2} in~\eqref{eqn:del2J_dot}, we obtain
\begin{align}
\frac{d}{d \tau}{\delta^{2}J}  = -\left( \frac{\partial^{2} B}{\partial \boldsymbol \zeta \partial \xi_{l}} \ddot{\xi}_{l}\right) \cdot \dot{\boldsymbol \zeta} \bigg|_{s=l}.
\end{align}
The R.H.S of this equation can be expressed as
\begin{align}
\label{eqn:Distinguished_biff_diagram_free_end}
 \begin{split}
 - \frac{\partial^{2} B}{\partial \boldsymbol \zeta \partial \xi_{l}} \ddot{\xi}_{l} \cdot \dot{\boldsymbol \zeta} &= - \ddot{\xi}_{l} \left(\frac{\partial}{\partial \boldsymbol \zeta}B_{\xi_{l}} \cdot \dot{\boldsymbol \zeta}  + \frac{\partial}{\partial \xi_{l}}B_{\xi_{l}} \dot{\xi} \right), \\ 
&=  -\ddot{\xi}_{l} \frac{d}{d \tau} B_{\xi_{l}}.
\end{split}
\end{align}
In this case, as a simple fold is traversed, the projection $B_{\xi_{l}}$ provides information on changes in the sign of critical eigenvalue. For a simple fold, $\ddot{\xi}_{l}$ and 
\begin{align}
\label{eqn:transversal_freeend}
\frac{d}{d \tau} B_{\xi_{l}}
\end{align}
 are non-zero. The later condition on ~\eqref{eqn:transversal_freeend} is equivalent to that given in~\eqref{eqn:transversal}. Consequently, the change in index is determined by the shape of  $B_{\xi_{l}}$ vs $\xi_{l}$ plot as illustrated in Figure~\ref{fig:schematic_bifurcation_diagram_b}. 

\section{Planar Elastica}
To test the theory developed in section~\ref{sec:Sec1}, we consider a special Euler elastica. In this section, we describe this apparatus and subsequently examine its straight equilibria and critical parameters.
\subsection{Setup}
\label{sec:Sec2}
\begin{figure}[t!]
    \centering
    \includegraphics[width=0.8\textwidth]{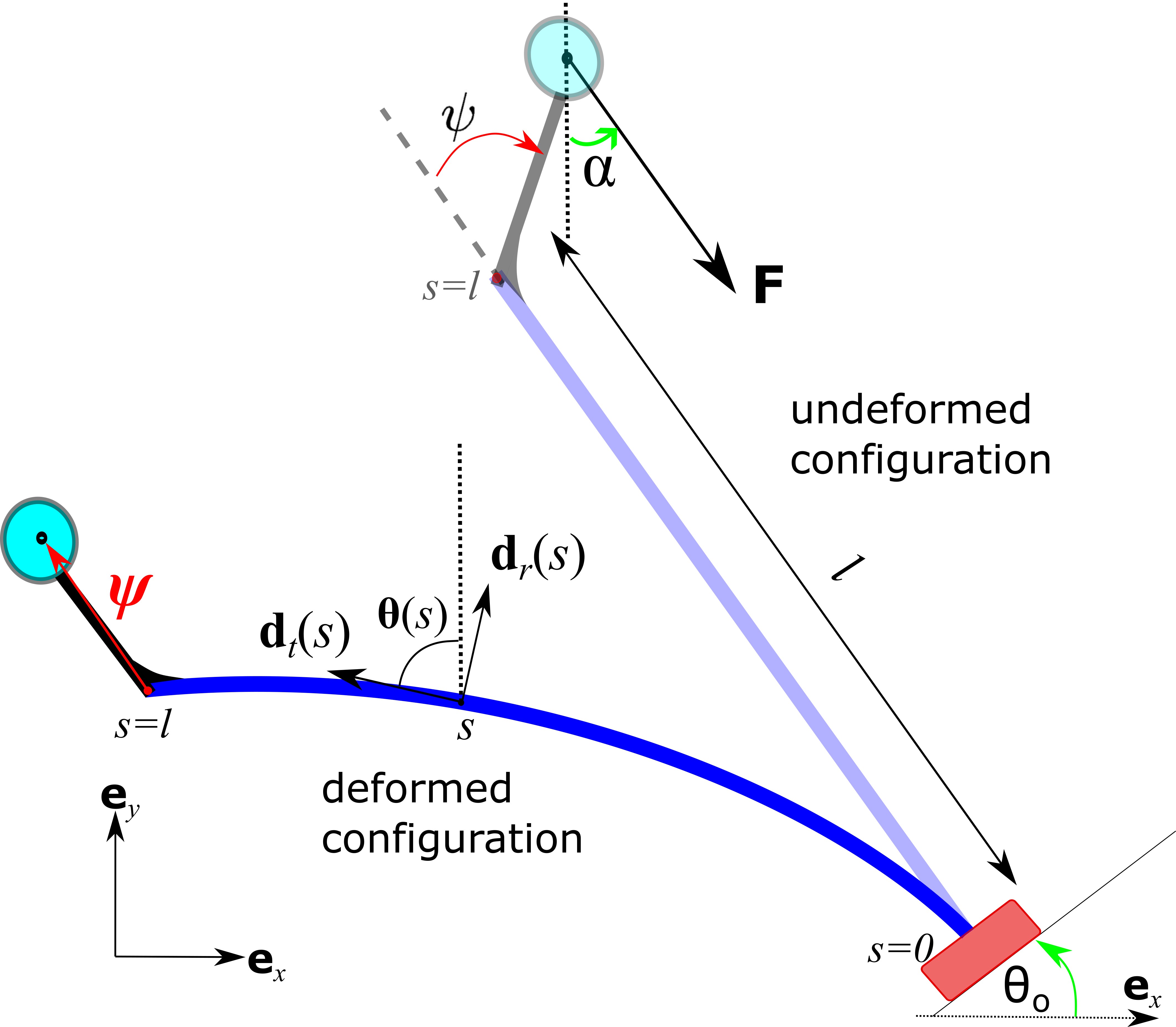}
    \caption{Schematic of an elastica with a load arm attached at the free end.}
    \label{fig:Elastica_schematic}
\end{figure}

An inextensible and unshearable planar massless uniform elastic rod with a bending stiffness $K$, which is straight in its undeformed configuration, is clamped at one end and connected to a dead load through a massless rigid lever arm at the other end, as shown in Figure~\ref{fig:Elastica_schematic}. Let $ \{ \mathbf{e}_{x},\mathbf{e}_{y}, \mathbf{e}_{z} \}$ denote a fixed right-handed orthonormal laboratory frame with $\mathbf{e}_{z}$ pointing out of the plane. The ends of the rod are defined by the coordinates $s=0$ and $s=l$ in its fixed undeformed configuration. The lever arm is rigidly fixed to the tip $s=l$ so that it forms an angle $\psi$ with the tangent of the tip $s=l$.  The dead load exerts a force $\mathbf{F}$ in the $\mathbf{e}_{x}-\mathbf{e}_{y}$ plane and a moment $ \boldsymbol \Psi \times \mathbf{F}$ due to an offset of the load $\mathbf{F}$ at $s=l$. Here, $\boldsymbol \Psi$ is the lever arm of the load. \\

The position vector of the rod centerline $\mathbf{r}(s)$ in the deformed configuration is represented using 
\begin{align}
\mathbf{r}(s)=  x(s) \mathbf{e}_{x}+ y(s) \mathbf{e}_{y}.
\end{align}
Let $\theta(s)$ be be the angle between the tangent of the rod
centerline and the vertical. Let $\mathbf{d}_{r}$ and $\mathbf{d}_{t}$ denote the local directors of the material cross-section and are represented by
\begin{align}
\mathbf{d}_{t}(s)=-\sin \theta(s) \mathbf{e}_{x} + \cos \theta(s)  \mathbf{e}_{y} , \qquad \mathbf{d}_{r}(s)=\cos \theta(s)  \mathbf{e}_{x} + \sin \theta(s)  \mathbf{e}_{y}.
\end{align}
The inextensibility and unshearability of the rod constrains the director $\mathbf{d}_{t}$ along the tangent of the centerline resulting in
\begin{align}
\label{eqn:Inextensibility_unshearability}
 x^{\prime}(s)&=-\sin \theta(s), \qquad y^{\prime}(s)=\cos \theta(s).
\end{align} 
Then, $\boldsymbol \Psi$ is given by $\Delta_{r}\mathbf{d}_{r}(l) + \Delta_{t}\mathbf{d}_{t}(l)$, where the constants $\Delta_{r}$ and $\Delta_{t}$ are the respective components. If $\Delta$ denotes the length of the arm and makes an angle $\psi$ with the tangent at the tip $s=l$, then $\Delta_{t}=\Delta \cos \psi$ and $\Delta_{r}=\Delta \sin \psi$. Using the principles of differential geometry, the curvature of the rod's centerline can be easily expressed as $\theta^{\prime}$. Let $\mathbf{m}(s)$ denote the internal moment vector along the rod. We consider a simple linear constitutive law where the planar bending moment in the rod is related to the curvature $\theta^{\prime}$ through
\begin{align}
\mathbf{m} \cdot \mathbf{e}_{z}=K \theta^{\prime},
\end{align}
    and the stored bending energy satisfies 
    \begin{align}
 \frac{1}{2}K \theta^{\prime ^{2}}.
    \end{align}
    The clamped end is rotated about $\mathbf{e_{z}}$- axis. Since it is a planar case, the moment vector $\mathbf{m}$ is always oriented perpendicular to the plane $(\mathbf{m}\cdot \mathbf{e}_{x}= \mathbf{m}\cdot \mathbf{e}_{y}=0 )$. The total energy stored in the system is the sum of the bending energy and the work done by the external force, satisfying the inextensibility and unshearability constraint \eqref{eqn:Inextensibility_unshearability} 
\begin{align}
\label{eqn:Gravitational_functional}
\mathcal{L}=\int_{0}^{l} \frac{K}{2} \theta^{\prime 2} + n_{x} (x^{\prime} + \sin \theta)  + n_{y} (y^{\prime} - \cos \theta) ds - \mathbf{F} \cdot( \mathbf{r}(l)+ \boldsymbol \Psi (l)),
\end{align}
where $n_{x}$ and $n_{y}$ are the corresponding Lagrange multipliers and denote the internal forces along $\mathbf{e}_{x}$ and $\mathbf{e}_{y}$ axes, respectively. We represent them using the vector $\mathbf{n}=[n_{x}, n_{y}]^{T}$. The equilibrium equations are obtained by applying the Euler-Lagrange equations~\eqref{eqn:Euler_Lagrange_bif_dist} on~\eqref{eqn:Gravitational_functional}
\begin{subequations}
\begin{align}
n_{x}^{\prime}&=0, \qquad n_{y}^{\prime}=0,\\
-K \theta^{\prime \prime}&+ n_{x} \cos \theta + n_{y} \sin \theta=0,
\label{eqn:Moment_balance}
\end{align}
\label{eqn:equilibria}
\end{subequations}
and natural boundary conditions~\eqref{eqn:Second_order_cond} at the free end $s=l$
\begin{subequations}
\begin{align}
    \label{eqn:N_BC_1}
    n_{x}(l) &-F_{x}=0, \quad n_{y}(l)-F_{y}=0,\\
     \label{eqn:N_BC_2}
       &K \theta^{\prime}(l) - \left(\boldsymbol \Psi \times \mathbf{F}\right) \cdot \mathbf{e}_{z}=0 ,
\end{align}
\end{subequations}
where $F_{x}$ and $F_{y}$ are the components of $\mathbf{F}$ along the $\mathbf{e}_{x}$ and $\mathbf{e}_{y}$ directions, respectively. In fact, the Euler-Lagrange equations given by~\eqref{eqn:equilibria} are equivalent to planar versions of the force and moment balance in elastic rods~\cite{antman2006nonlinear}

\begin{align}
\begin{split}
\mathbf{n}^{\prime}&=\mathbf{0},\\
 \mathbf{m}^{\prime} + \mathbf{r}^{\prime} \times \mathbf{n}&=\mathbf{0}.
 \end{split}
 \label{eqn:Newton_balance}
\end{align}
The load vector $\mathbf{F}$ makes an angle $\alpha$ with respect to vertical as shown in Figure~\ref{fig:Elastica_schematic}. The load $\mathbf{F}$ can be represented using its magnitude and direction $\alpha$ as
\begin{align}
\mathbf{F}= | \mathbf{F} | \left( \sin \alpha \mathbf{e}_{x} - \cos \alpha \mathbf{e}_{y} \right).
\end{align}
 In this apparatus, the parameters that describe the Dirichlet boundary conditions, such as $\theta_{o}$ and those that characterize the natural boundary conditions such as $| \mathbf{F} |$, $\Delta$, $\alpha$ and $\psi$ can be varied. We perform several numerical experiments using parameter continuation \cite{DOEDEL1991a} and employ the proposed distinguished bifurcation diagrams to assess the stability of the equilibria.

\subsection{Bifurcation Analysis}
\label{ssec:Bifurcation}
The distinguished bifurcation diagram provides information on changes in stability. So, the stability of atleast one equilibrium along a branch must be established before proceeding. For this reason, we conduct bifurcation analysis on the current apparatus and establish stability indices for some equilibria. The buckling characteristics of the elastica are well known; now, let us examine how the lever arm influences this behavior. For this analysis, we set $\alpha=0$ and $\psi=0$, and consequently we have $n_{x}=0$ and $n_{y}=-| \mathbf{F} |$. The classical elastica is obtained when $\Delta=0$. Before proceeding further, we non-dimensionalize~\eqref{eqn:Moment_balance} by substituting $\bar{s}=\frac{s}{l}$, $\epsilon= \frac{\Delta}{l}$ and $P=\frac{\\|\mathbf{F}|l^{2}}{K}$ that fetches 
\begin{subequations}
\begin{align}
-\theta^{\prime \prime} & - P \sin \theta=0.\\
\theta(0)&=0, \qquad \theta^{\prime}(1) -  P \epsilon \sin \theta =0.
\end{align}
\end{subequations}
It has the trivial solution $\theta(s)=0$. On linearizing the equilibria about this solution, we get
\begin{subequations}
\begin{align}
\mathcal{S}\delta \theta &\equiv -\delta \theta^{\prime \prime} - P \delta \theta=0,\\
\delta \theta(0) &=0,  \qquad 
\delta \theta^{\prime}(1) - P \epsilon \delta \theta(1)=0 
\end{align}
\label{eqn:Linearized_Elastica}
\end{subequations}
where $\delta \theta(s)$ are perturbation in $\theta$ satisfying the boundary conditions. Its non-trivial solutions are of the form:
\begin{align}
\delta \theta(s)= A \sin \sqrt{P}s,
\end{align}
and the boundary condition at $s=l$ leads to the relation
\begin{align}
\begin{split}
&\sqrt{P} \cos \sqrt{P} - P \epsilon \sin \sqrt{P}=0, \\ 
&\implies \cot \sqrt{P} - \sqrt{P}\epsilon=0, \qquad  \sqrt{P} \neq n \pi.
\end{split}
\end{align}

\begin{figure}[t!]
    \centering
    \includegraphics[width=1.0\textwidth]{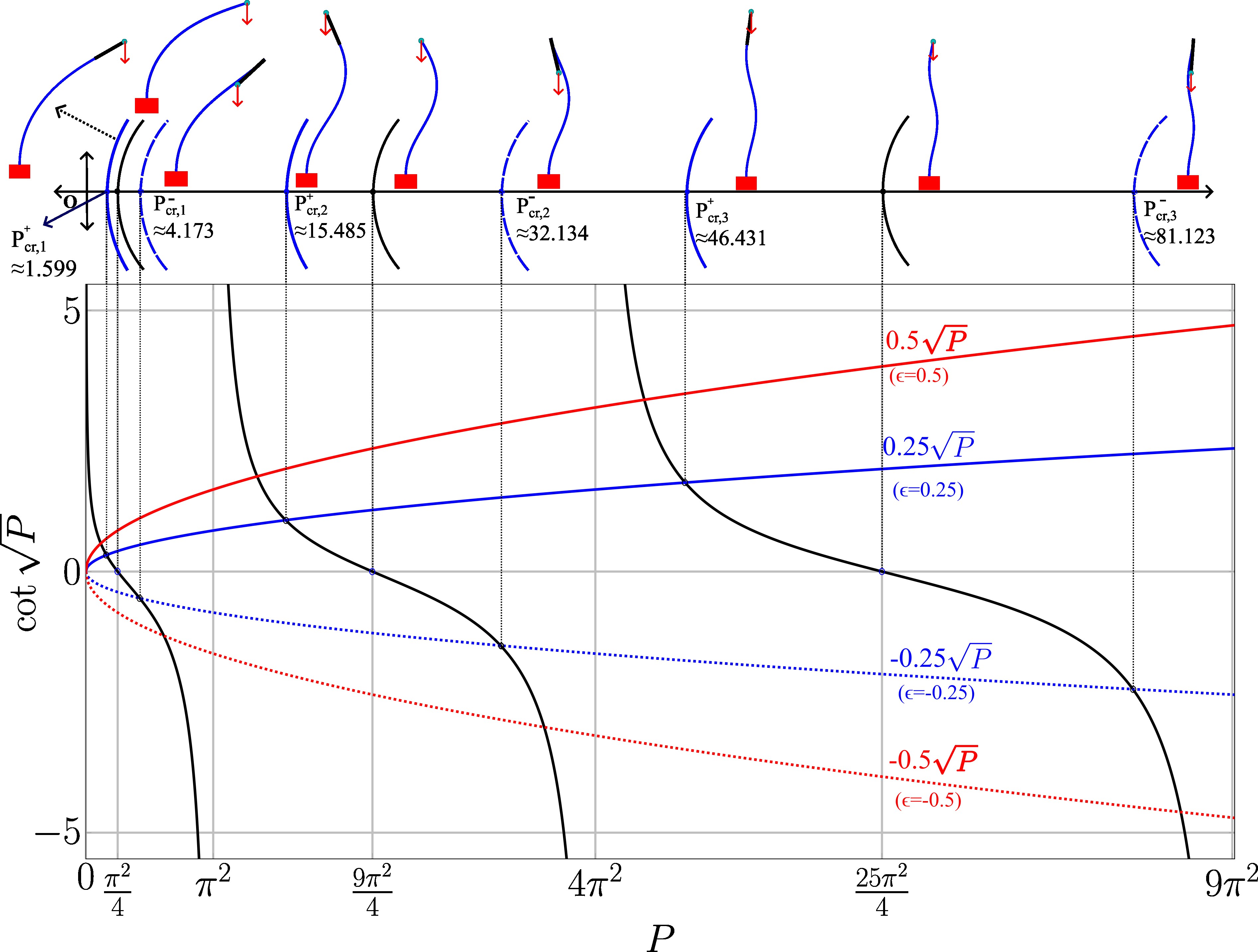}
    \caption{Plot illustrating the solutions of $\cot{\sqrt{P}} - \sqrt{P}\epsilon=0$ for different values of $\epsilon$, along with the local bifurcation characteristics associated with $\epsilon=+/-0.25$ shown at top. The family of configurations branching out at bifurcations are indicated.}
    \label{fig:schematic_bifurcation_diagram}
\end{figure}

The non-trivial solutions satisfying this relation that correspond to bifurcations are depicted in Figure~\ref{fig:schematic_bifurcation_diagram}. The critical buckling modes of elastica are the solutions pertaining to $\epsilon=0$ denoted as $P^{o}_{cr}=\frac{(2k-1)^{2}\pi^{2}}{4}, k=1,2,\dots,\infty$. Let $P_{cr}^{+}$ represent the critical load values pertaining to positive $\epsilon$, while $P_{cr}^{-}$ represent pertaining to negative $\epsilon$. Note that for each $k$ th mode, $ (k-1)^{2} \pi^{2}<P_{cr,k}^{+}< \left( \frac{2k-1}{2}\right)^{2} \pi^{2}$ and $\left( \frac{2k-1}{2}\right)^{2} \pi^{2}<P_{cr,k}^{-}< k^{2} \pi^{2}$, where $k=1,2,\dots,\infty$. It can be easily shown that these points correspond to supercritical pitchfork bifurcations, following a similar analysis in ~\cite[pp,~174-177]{antman2006nonlinear}. The linearized equilibrium~\eqref{eqn:Linearized_Elastica} can be rewritten as an eigenvalue problem: 
\begin{align}
 -\delta \theta^{\prime \prime} - P \delta \theta= \mu \delta \theta , \qquad \delta \theta(0) =0,  \qquad 
\delta \theta^{\prime}(1) - P \epsilon \delta \theta(1)= 0.
\end{align}
where $\mu$ is the eigenvalue. On rearranging the terms we have
\begin{align}
 \delta \theta^{\prime \prime} + \left(P + \mu \right)\delta \theta= 0,\qquad \delta \theta(0) =0, \qquad 
\delta \theta^{\prime}(1) -  P  \epsilon \delta \theta(1)= 0.
\end{align}
Typically, at a supercritical bifurcation~\cite{iooss2012elementary}, the number of negative eigenvalues increases by one with each successive bifurcation. In our setup, the index of the straight equilibria (trivial solutions) for load $P_{cr,k}^{+}< P < P_{cr,k+1}^{+}$  for positive $\epsilon$ (and $P_{cr,k}^{-}< P < P_{cr,k+1}^{-}$ for negative $\epsilon$) is $k$. 
For a very high $\epsilon$ value, i.e., when the $\sqrt{P}\epsilon$ line approaches the vertical line, the value of $P_{cr}^{+}$ approaches zero, suggesting that bifurcations occur at a very small load when the arm is long. On the other hand, $P_{cr}^{-}$ approaches $\pi^{2}$ for higher loads, which indicates the critical loads for the case of fixed-fixed ends. The trivial solutions are stable only for loads below $P_{cr,1}^{+}$ for positive $\epsilon$ (and below $P_{cr,1}^{-}$ for negative $\epsilon$). Consider two instances of loads which correspond to first mode and second mode, say $\pi^{2}/4$ and $9 \pi^{2}/4$. These chosen straight equilibria correspond to index one and two respectively, while the negative $\epsilon$ leads to index zero and index one respectively. We use these equilibria as initial solutions for numerical continuation in the subsequent sections.

\section{Numerical Examples} 
\label{sec:Sec3}
The distinguished bifurcation diagrams facilitate the determination of the index change at a fold in the varying parameter. If there are no folds, the stability remains the same for all solutions along the branch. At folds, the stability transitions usually occur, and they were evaluated by examining the manner in which the fold occurs, as depicted in Figure~\ref{fig:Distinguished_Bifurcation}. We use the straight trivial solutions from the previous section as initial solutions to generate the family of equilibria by varying parameters. Subsequently, we assess the stability of the resulting equilibria. Numerical continuation is performed using AUTO-07p~\cite{doedel2007auto}, which employs pseudo-arclength continuation~\cite{DOEDEL1991a}, and is capable of detecting singularities like folds and bifurcations as the parameter is varied. In these examples, the analysis is restricted solely to simple folds.

\subsection{Varying Parameter in Free end}
\label{ssec:Examples_free_end}
In the first subsection, the parameters associated with the free end $(s=l)$ are varied. These include the magnitude of the tip load $P$, direction of the tip load $\alpha$, length of the load arm $\epsilon$, and the orientation of the arm $\psi$. Moreover, the parameters $P$ and $\epsilon$ can assume negative values, and they indicate the negative direction of the force vector $\mathbf{F}$ and the arm vector $\boldsymbol \Psi$ respectively. In the presented examples, stiffness $K$ and length $l$ are set to unit values and, therefore $P = |\mathbf{F}|$ and $\epsilon = \Delta$.

\subsubsection{Rotating Arm}
\label{sssec:Example1}
We begin with a numerical example in which continuation is performed with respect to $\psi$, simulating the quasi-static rotation of the arm. The system is $2\pi$-periodic in $\psi $, i.e., the system corresponding to $\psi =k$ and $\psi=2 \pi + k$ exhibit identical features for any real $k$. Consequently, the solution family also exhibits $2 \pi$- periodic behavior. The analysis begins with parameters $P=\pi^{2}/4 $ and $\epsilon=0.25$, and continuation is initiated from the straight trivial equilibrium at $\psi=0$ which has an index of one ($P_{cr,1}^{+}< P < P_{cr,2}^{+} $ in Section~\ref{ssec:Bifurcation}). When $\psi$ is the varying parameter, the ordinate in the distinguished bifurcation diagrams, as specified by~\eqref{eqn:Distinguished_biff_diagram_free_end} is
\begin{align*}
\frac{\partial B}{\partial \psi} &= -\frac{\partial}{\partial \psi} \left(\mathbf{F} \cdot (\mathbf{r}(l) + \epsilon \cos \psi \mathbf{d}_{t}(l) + \epsilon \sin \psi \mathbf{d}_{r}(l) ) \right),\\
& = \left(\boldsymbol \Psi(l) \times \mathbf{F} \right) \cdot \mathbf{e}_{z} = \mathbf{m}(l) \cdot \mathbf{e}_{z}.
\end{align*}

\begin{figure}[t!]
     \centering
     \begin{subfigure}{0.9\textwidth}
         \centering
         \includegraphics[width=1.0\textwidth]{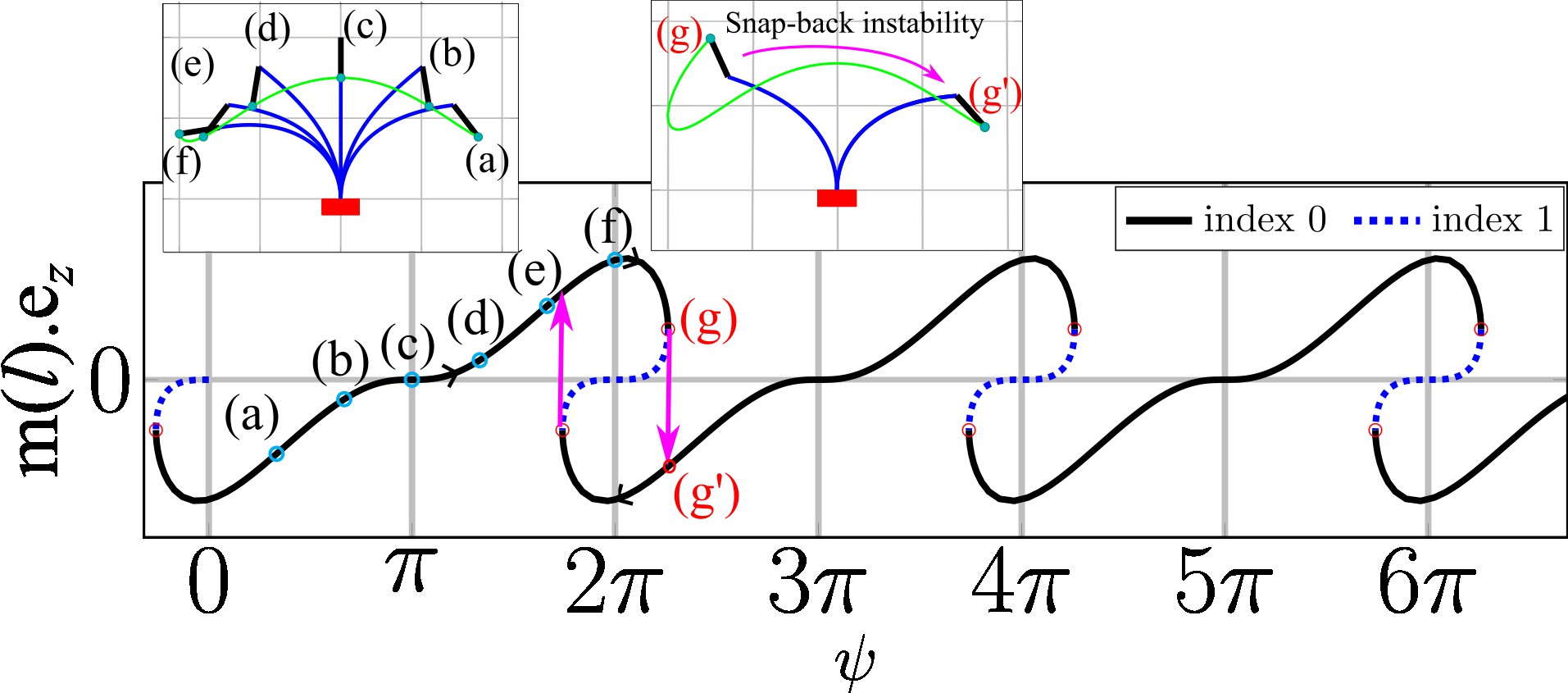}
         \caption{$P=\frac{\pi^{2}}{4}$}
         \label{fig:Example_1_a}
     \end{subfigure}
    \begin{subfigure}{0.9\textwidth}
         \centering
         \includegraphics[width=1.0\textwidth]{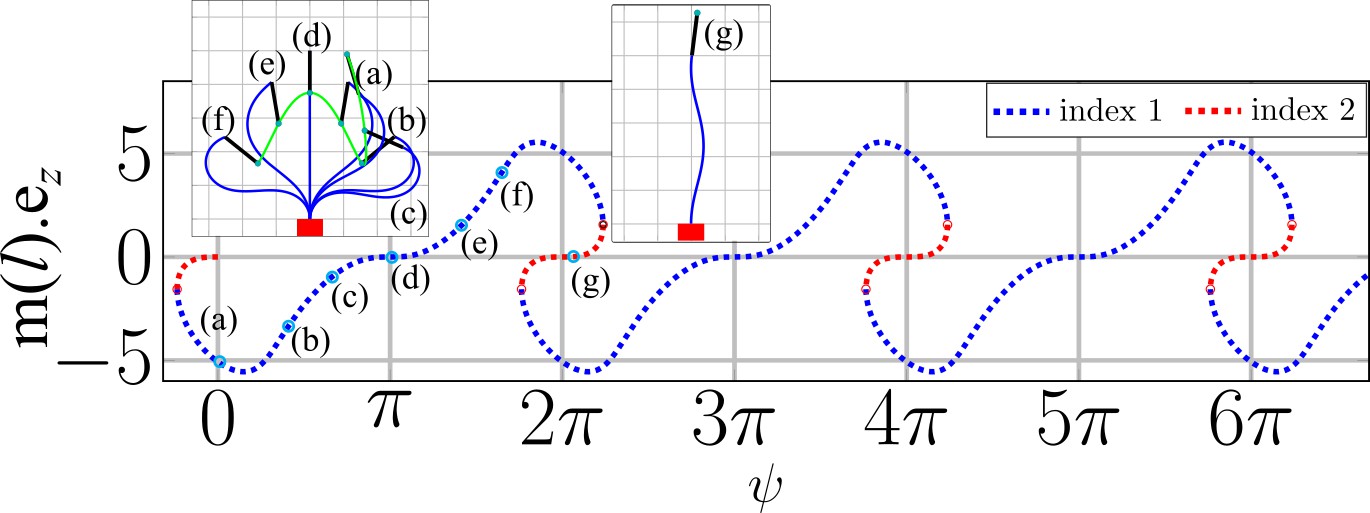}   
         \caption{$P=\frac{9\pi^{2}}{4}$}
         \label{fig:Example_1_b}
     \end{subfigure} 
     \caption{The distinguished bifurcation diagram for arm $\epsilon = 0.25$ rotating for loads (a)$P=\frac{\pi^{2}}{4}$ and (b) $P=\frac{9\pi^{2}}{4}$. The intermediate configurations along the family of stable equilibrium are labeled and depicted at the top. The stable configurations before and after the snap-back instability are also displayed. In (b), all equilibria are unstable. Although the folds are present no information on snapping behavior can be inferred, as no stable equilibrium exists. The tip trace in the foldless region is shown (in green).}
     \label{fig:Example_1}
\end{figure}
In the last step, the boundary condition~\eqref{eqn:N_BC_2} is applied. This ordinate $\mathbf{m}(l).\mathbf{e}_{z}$ is evaluated from the numerical continuation solutions and plotted against the parameter $\psi$ to generate the associated distinguished bifurcation diagram, as shown in Figure~\ref{fig:Example_1}. This plot exhibits folds, indicating stability transitions. The direction of change in index near folds is inferred from Figure~\ref{fig:schematic_bifurcation_diagram_b} and is indicated in Figures~\ref{fig:Example_1}. The equilibrium at $\psi=\pi$ is equivalent to the case of $\epsilon<0$ in Section~\ref{ssec:Bifurcation}, where $ P< P_{cr,1}^{-}$ corresponds to index zero and agrees with the diagrams. The dynamical aspects of the snapping motion are not discussed in this paper, and for further details, refer to~\cite{snyder1990dynamics, Armanini2017}. The solutions and their projections exhibit $2 \pi$- periodic characteristics. Figure~\ref{fig:Example_1_a} also displays the elastica configurations at intermediate values of $\psi$. The potential snapping motion depicting the initial and final configurations is also shown. This motion depends on the current state and the direction of current parameter change, and the snap jump is indicated by a red arrow.
 
 Figure~\ref{fig:Example_1_b} displays the bifurcation diagrams for a family of equilibria at a higher load, $P=9 \pi^{2}/4$. Under this loading, the straight trivial equilibrium has an index of two ($P_{cr,2}^{+}< P < P_{cr,3}^{+} $ in Section~\ref{ssec:Bifurcation}). These plots also exhibit folds, which, although appear qualitatively similar to the case of $P=\pi^{2}/4$ but differ significantly in their stability characteristics. In this scenario, only families of equilibria with indices one and two exist, corresponding to unstable equilibria. To obtain a family of stable equilibria, the continuation process must be initiated from a different solution.  The potential snap-back instability can be assessed only when a family of stable equilibria exists.

\begin{figure}
     \centering
     \begin{subfigure}{0.9\textwidth}
         \centering
         \includegraphics[width=1.0\textwidth]{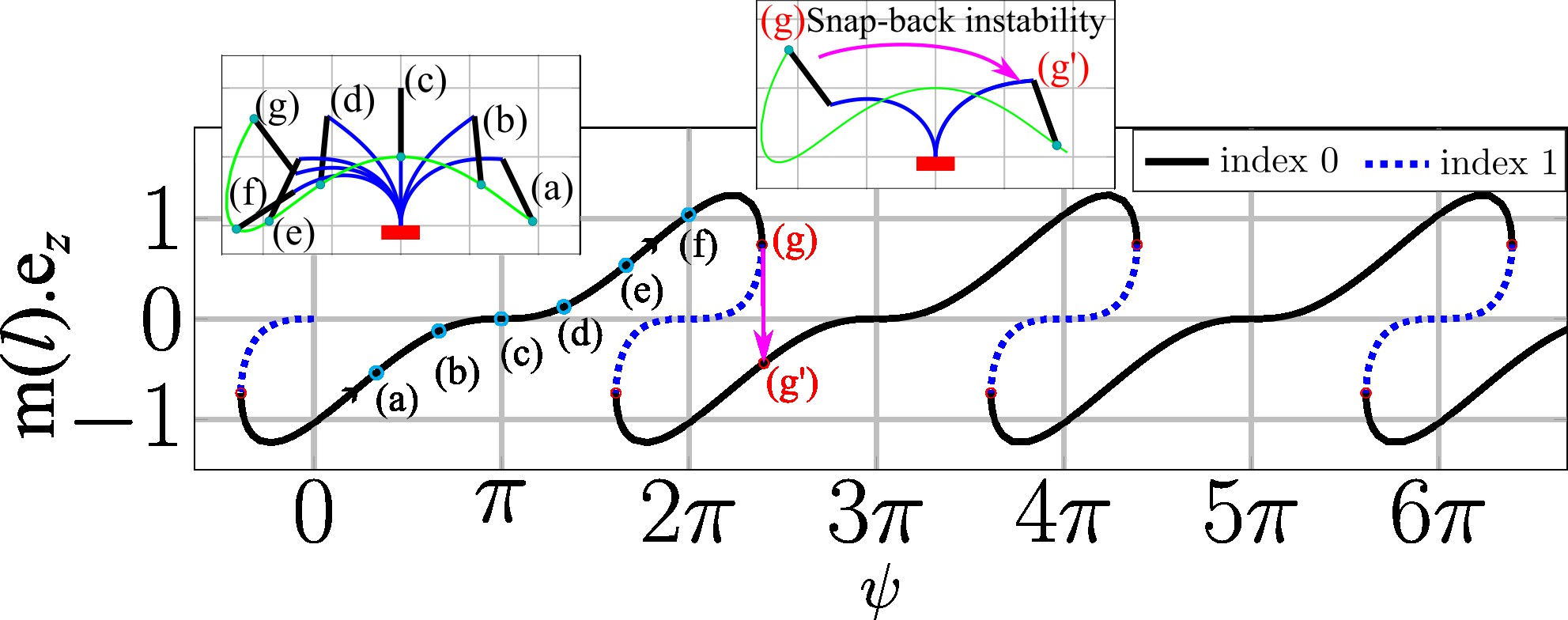}
         \caption{$P=\frac{\pi^{2}}{4}$}
         \label{fig:Example_1_c}
     \end{subfigure}
    \begin{subfigure}{0.9\textwidth}
         \centering
         \includegraphics[width=1.0\textwidth]{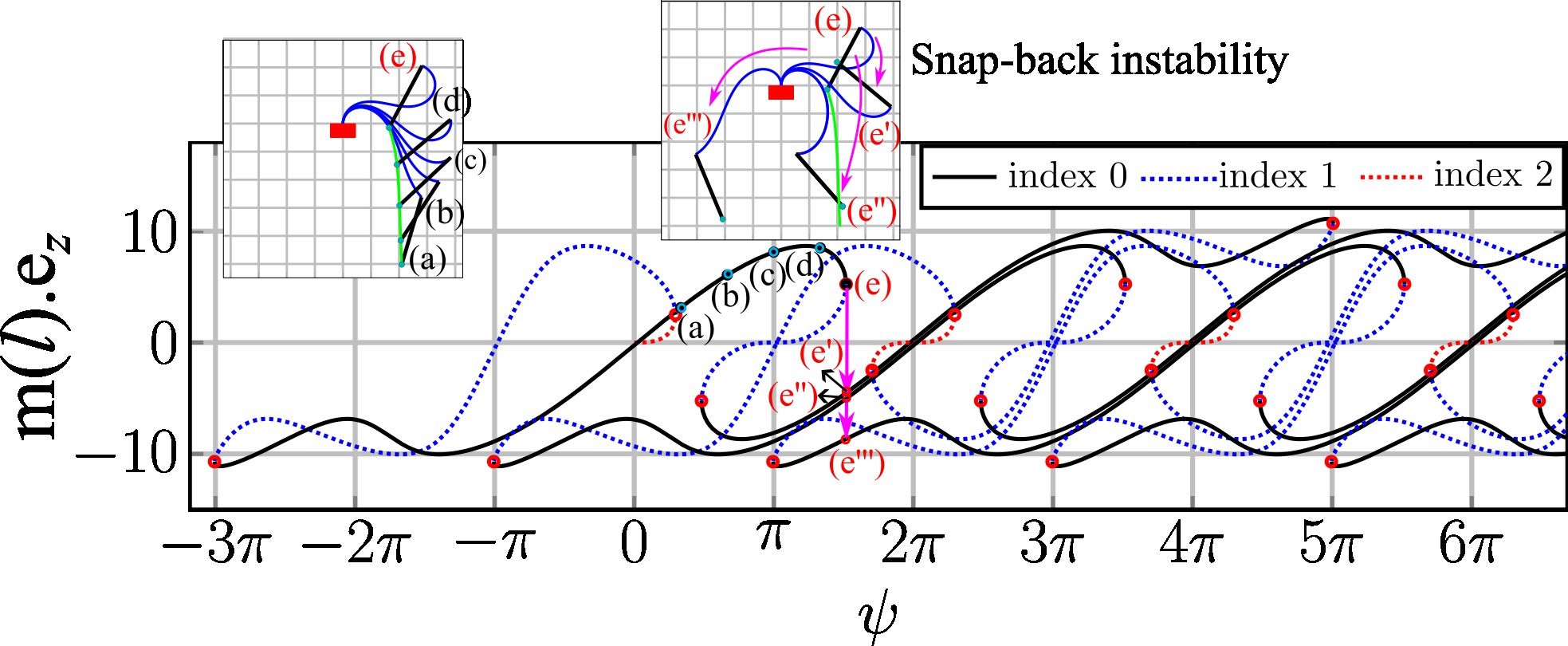}
         \caption{$P=\frac{9\pi^{2}}{4}$}
         \label{fig:Example_1_d}
     \end{subfigure}
          \caption{The distinguished bifurcation diagram for an arm $\epsilon = 0.5$ rotating for loads (a) $P=\frac{\pi^{2}}{4}$ and (b)$P=\frac{9\pi^{2}}{4}$. Intermediate configurations along this family are labeled and depicted at the top. Additionally, the stable configurations before and after the snap-back instability are displayed. The trajectory of the arm's tip preceding the snap-back instability is indicated in green.}
          \label{fig:Example_1_2}
\end{figure}

\begin{figure}[t!]
     \centering
     \begin{subfigure}{0.75\textwidth}
         \centering
         \includegraphics[width=1.0\textwidth]{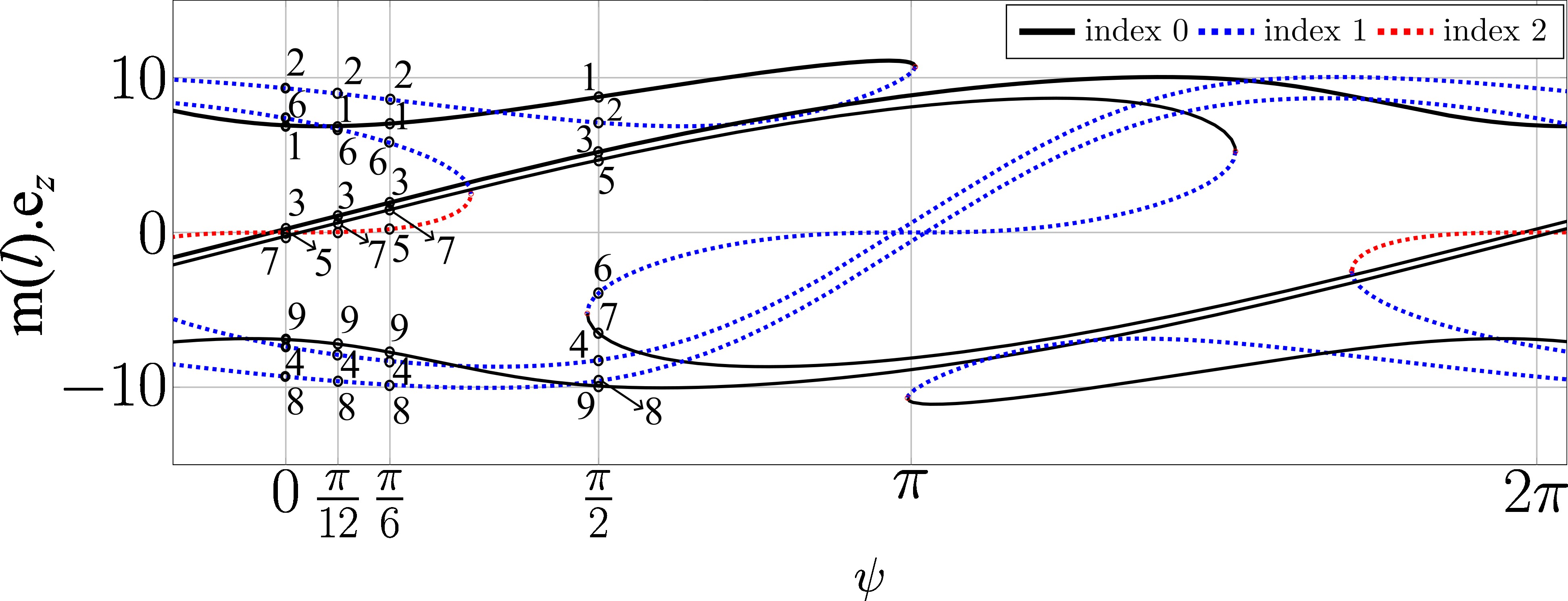}
         \caption{}
         \label{fig:example2_screenshot}
     \end{subfigure}
          \begin{subfigure}{0.9\textwidth}
         \centering
         \includegraphics[width=1.0\textwidth]{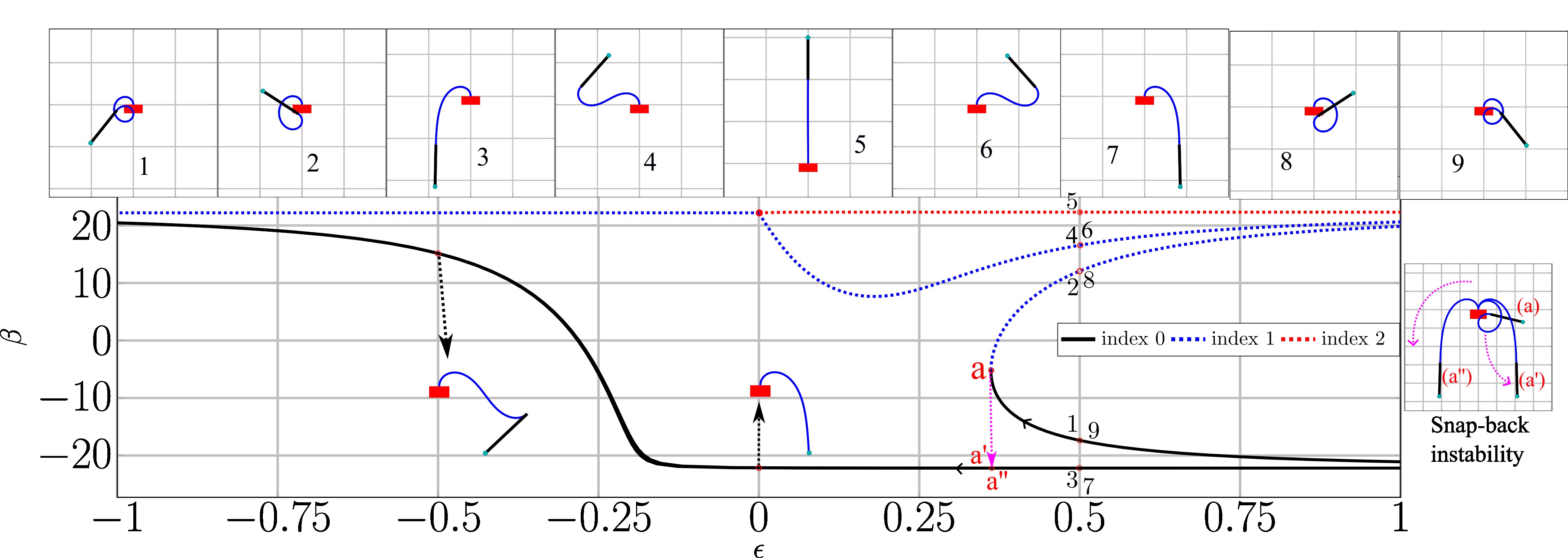}
         \caption{$ \psi=0$}
         \label{fig:Example2_0}
     \end{subfigure}
             \caption{(a) The bifurcation diagram from the previous analysis adjusted between $0$ and $2 \pi$. Equilibria are labeled according to the order which continuation traverses the specified $\psi$. (b) Distinguished bifurcation diagrams for $\epsilon$ varying from $0.5$ to $-0.5$ at $\psi=0$, showing the presence of folds and bifurcations. The equilibria corresponding to the given labels are shown at the top. The potential snap-back instability is also illustrated.}
        \label{fig:Example2_a}
\end{figure}

A similar analysis is performed for a longer arm $\epsilon =0.5$, and Figure~\ref{fig:Example_1_2} displays the corresponding plots. The straight trivial solutions with this arm have the same index as those of $\epsilon =0.25$ i.e, one and two for loads $P=\pi^{2}/4$ and $P=9\pi^{2}/4$, respectively. The bifurcation diagram and the configurations for $P=\pi^{2}/4$, shown in Figure~\ref{fig:Example_1_c} appear qualitatively similar to those in the case of $\epsilon=0.25$ displayed in Figure~\ref{fig:Example_1_a}. However, the higher load case of $P=9\pi^{2}/4$ exhibits significant differences. Here, we observe several folds, with more than two equilibria existing for an identical value of $\psi$ suggesting the presence of multi-stability. Unlike the previous example for $P=9 \pi^{2}/4$ and $\epsilon =0.25$ in Figure~\ref{fig:Example_1_b}, a family of equilibria with index zero exists here. Consequently, snap-back instability occurs as $\psi$ is varied beyond the folds. The plots also illustrate several equilibrium configurations during this maneuver. Furthermore, the stable equilibria before and after the set of snap-back instability are also displayed. We assume that the geometry of the lever arm and elastica allows for self-intersection, and the resulting configurations as displayed in the plots are valid. For a zero arm, when the load $P$ exceeds the critical buckling length, it buckles into a non-straight configuration, deflecting  either to the right or to the left, with both configurations being mirror images of each other. The non-zero asymmetrical arm also results in similar post-buckling structures, but they are not mirror images. The arm acts as an imperfect parameter breaking the symmetry. Furthermore, the arm also induces bi-stability in the buckled configurations, meaning there is more than one stable equilibria for some values of $\psi$ in buckled states, as shown in Figure~\ref{fig:Example_1_d}. Subsequently, more than two stable equilibria exist in some instances. When $\psi$ is varied past the fold, snap-back instability arises and leads to one of the available stable equilibria. The specific equilibrium reached depends on the amount of energy released during snapping and the dissipation of the system and is beyond the scope of this paper. With appropriate technical considerations, novel snapping, triggering, and switching mechanisms can be realized. The periodicity of the system with respect to $\psi $ allows the reproducibility of the snap-back instability, i.e., the similar snapping motion can be replicated by tuning the parameter forward by $ 2 \pi $. Furthermore, snapping motion can be generated by reversing the parameter $\psi$, but occurs at a fold in the reverse direction. 

\begin{figure}[t!]
     \centering
     \begin{subfigure}{0.9\textwidth}
         \centering
         \includegraphics[width=1.0\textwidth]{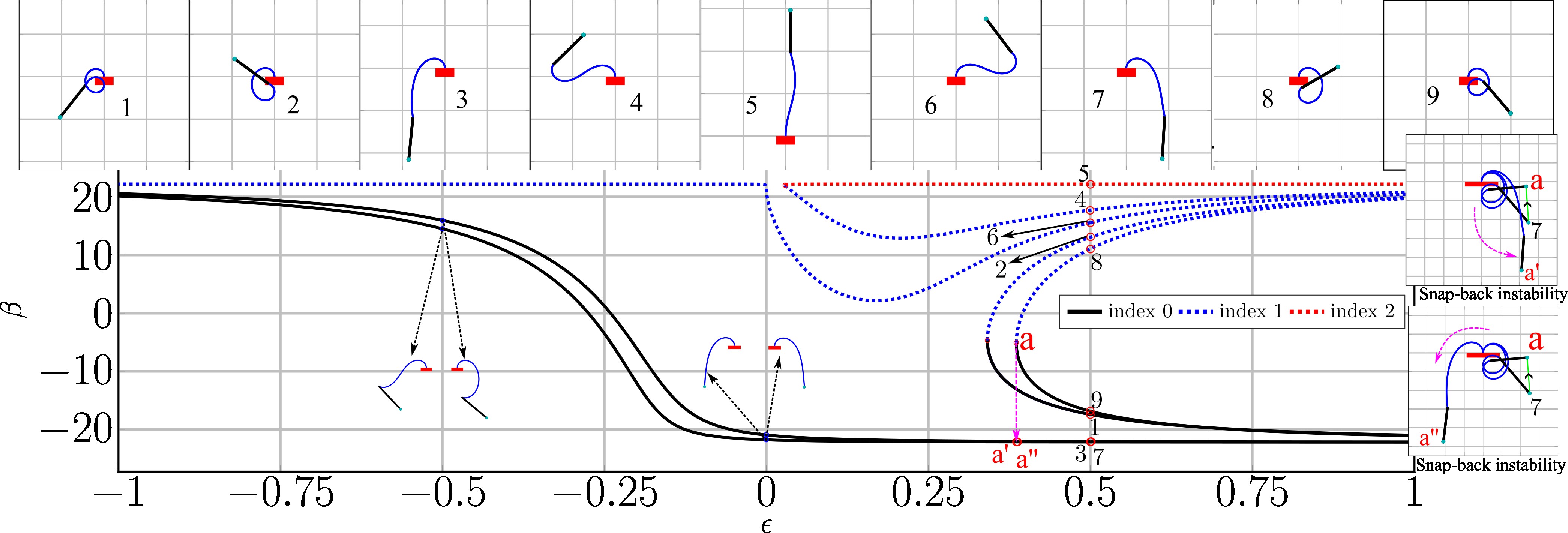}
         \caption{$\psi=\frac{\pi}{12}$}
         \label{fig:Example2_b}
     \end{subfigure}
          \begin{subfigure}{0.9\textwidth}
         \centering
         \includegraphics[width=1.0\textwidth]{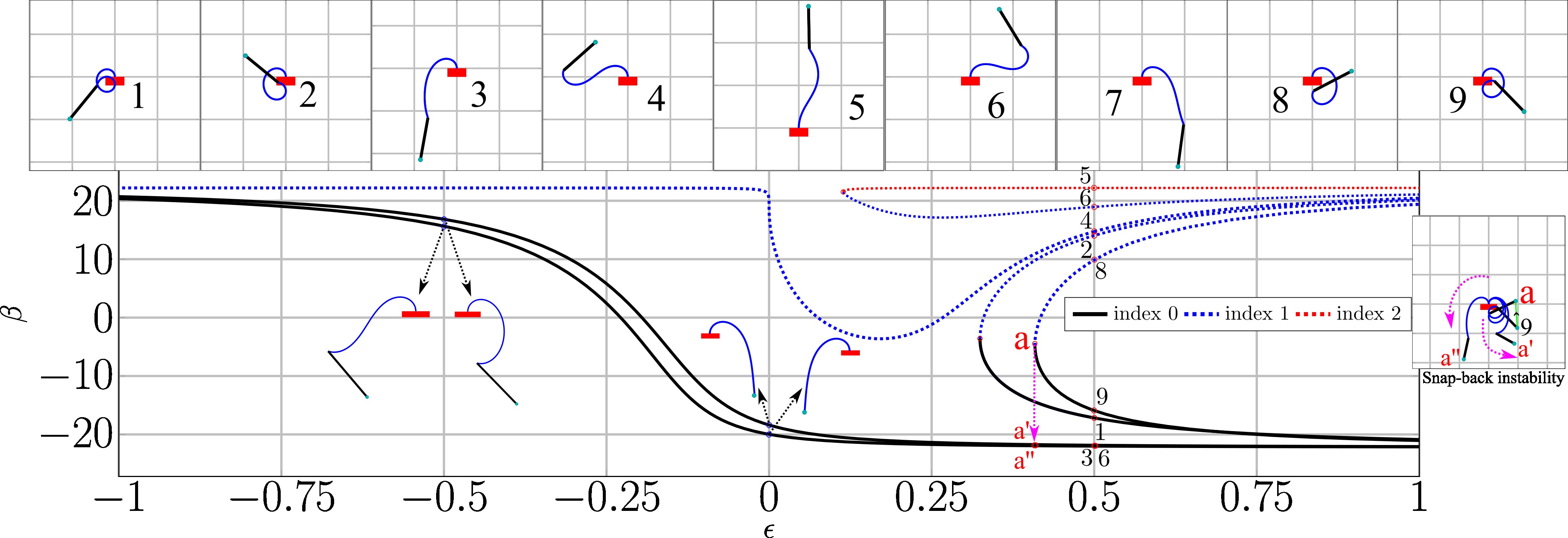}
         \caption{$ \psi=\frac{\pi}{6}$}
         \label{fig:Example2_c}
     \end{subfigure}

     \begin{subfigure}{0.9\textwidth}
         \centering
         \includegraphics[width=1.0\textwidth]{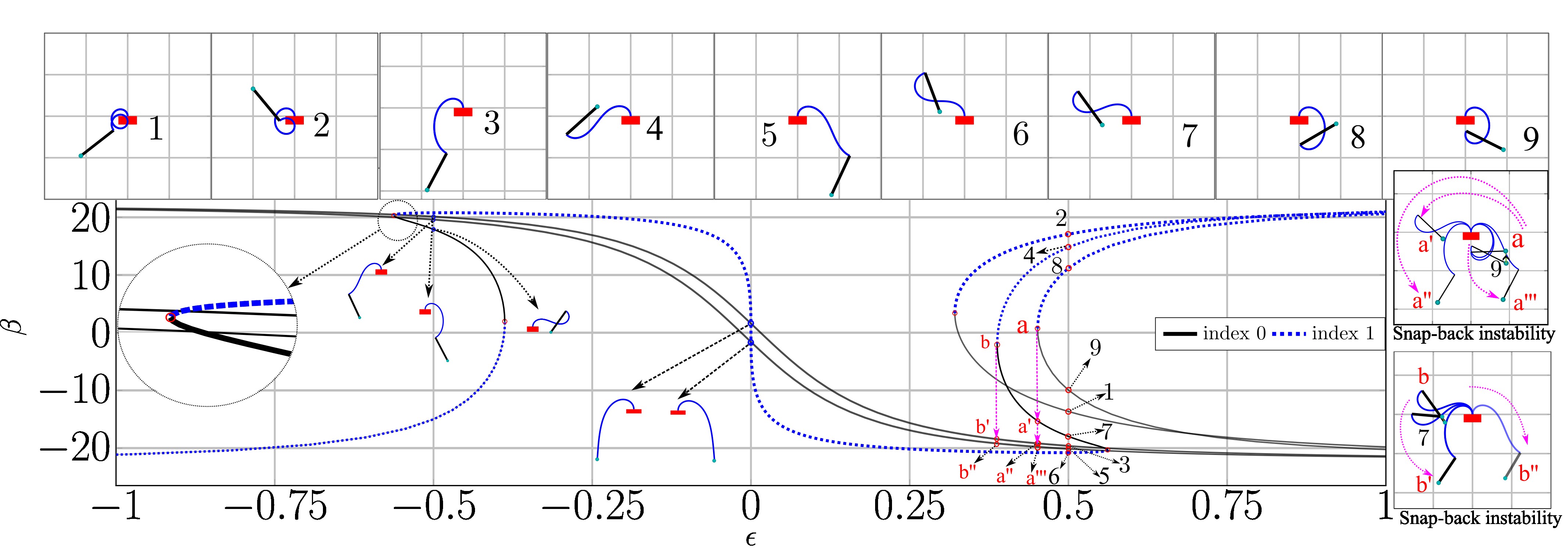}
         \caption{$\psi=\frac{\pi}{2}$}
         \label{fig:Example2_d}
     \end{subfigure}
        \caption{Distinguished bifurcation diagrams for the arm length $\epsilon$ varying from from $0.5$ to $-0.5$ at arm angles (a) $\psi=\pi/12$, (b) $\psi=\pi/6$, and (c) $\psi=\pi/2$. Several instances of snap-back instability can be noticed and a few of them are illustrated. An enlarged view of a section of plot (c) is provided to clearly illustrate the folds.}
        \label{fig:Example2_2}
\end{figure}

\begin{figure}[t!]
     \centering
     \begin{subfigure}{0.9\textwidth}
         \centering
         \includegraphics[width=\textwidth]{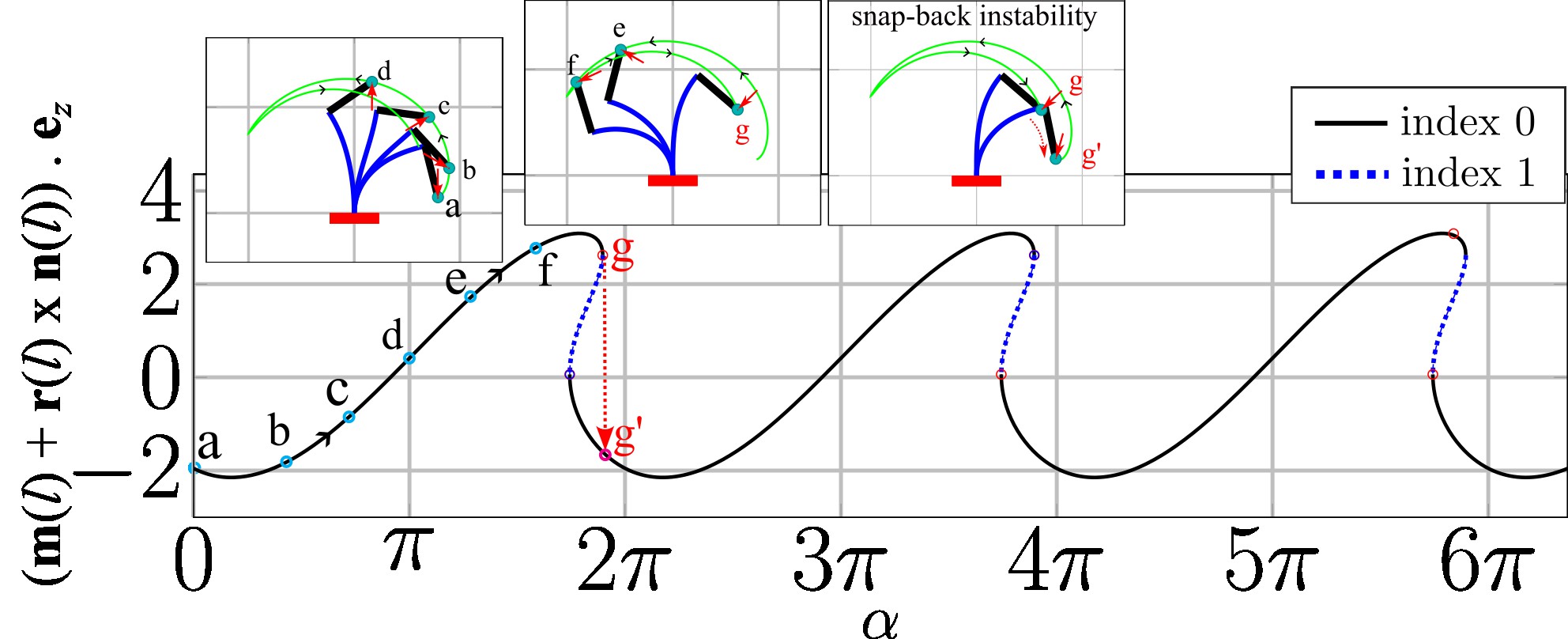}
         \caption{$P=\frac{\pi^{2}}{4}$}
         \label{fig:Example3c}
     \end{subfigure}
          \begin{subfigure}{0.9\textwidth}
         \centering
         \includegraphics[width=\textwidth]{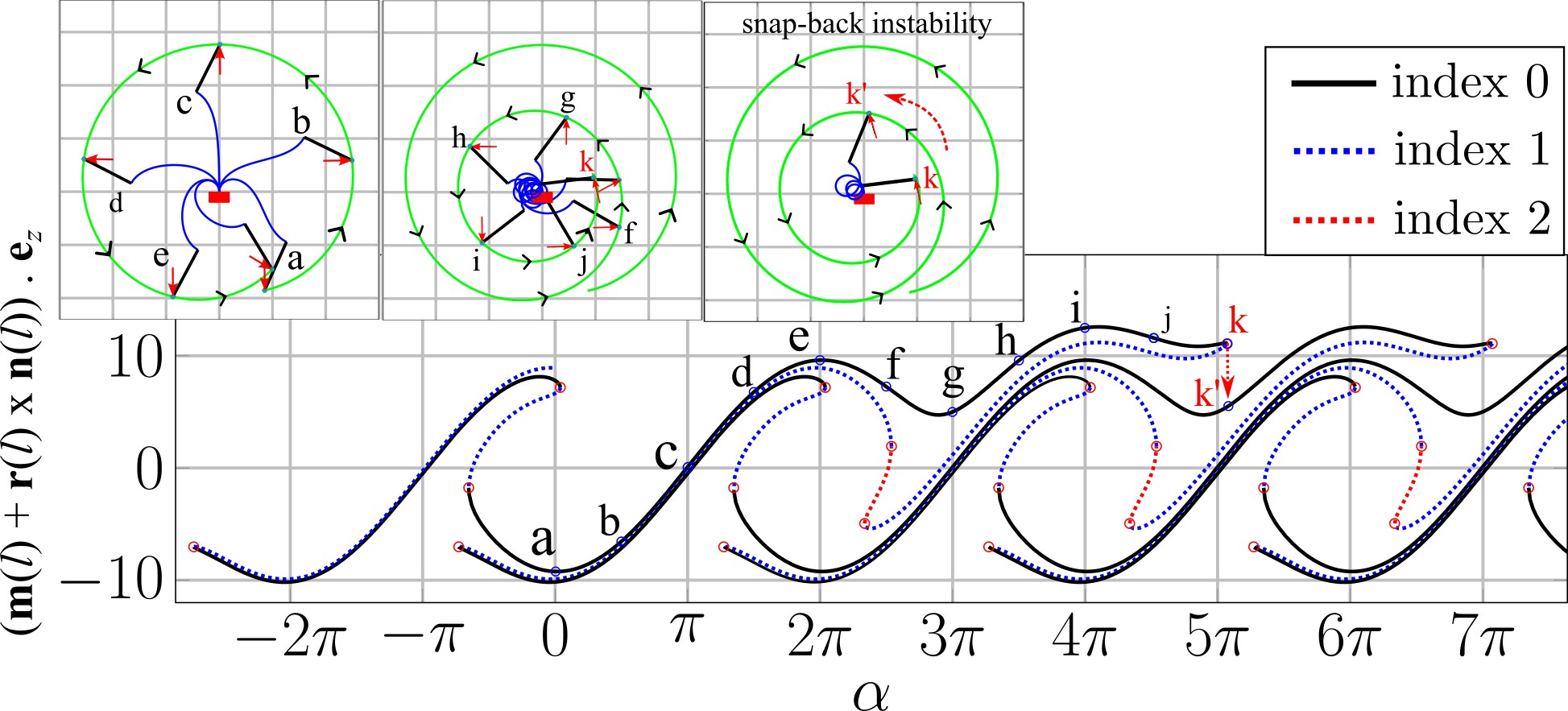}
         \caption{$P=\frac{9\pi^{2}}{4}$}   
         \label{fig:Example3d}
     \end{subfigure}
        \caption{Distinguished Bifurcation Diagrams when load $P$ with an arm $\epsilon=0.5$ varies its direction $\alpha$, completing several full rotations for $P=\frac{\pi^{2}}{4}$ and $P=\frac{9\pi^{2}}{4}$. The intermediate configurations along this stable-family are displayed at the top. The tip trace before the snap-back instability is also displayed (in green).}
        \label{fig:Example3}
\end{figure}

 \subsubsection{Varying Arm Length}
 \label{sssec:Example2}
We continue this analysis by varying the parameter associated with the arm's length $\epsilon$. The equilibria with $\epsilon=0.5$ at $\psi=0$, $\psi=\pi/12$, $\psi=\pi/6$, and $\psi=\pi/2$ from the previous analysis are chosen as initial solutions, and continuation is performed along $\epsilon$ from $0.5$ to $-0.5$. Figure~\ref{fig:example2_screenshot} displays the bifurcation plot from the previous analysis, truncated between $0$ and $2\pi$ to obtain all possible equilibria at a given $\psi$. Owing to the $2 \pi$-periodicity of the system, the section of the plot between $4 \pi$ and $6\pi$ is represented within the interval $0$ to $2 \pi$.  In this case, the ordinate of the distinguished bifurcation diagram is
\begin{align}
\beta:= \frac{\partial B}{\partial \epsilon}= -\mathbf{F} \cdot (\cos \psi \mathbf{d}_{t}  + \sin \psi \mathbf{d}_{r})= -\mathbf{n}(l) \cdot \mathbf{d}_{t}(l) \cos \psi  - \mathbf{n}(l) \cdot \mathbf{d}_{r}(l) \sin \psi,
\end{align}
which is the force component along the arm vector. Figure~\ref{fig:Example2_0} and Figure~\ref{fig:Example2_2} display the $\beta$ vs. $\epsilon$ plots associated with various $\psi$. A few plots exhibit folds and stability information is encoded according to Figure~\ref{fig:schematic_bifurcation_diagram_b}. Let us first consider the symmetric case of $\psi=0$. Here, the ordinate $\beta$, the load vector's projection along the arm $\Psi$, has identical values for deflections to both the right and left, resulting in overlapping plots. Consequently, the supercritical pitchfork bifurcation diagram appears only as a half plot, which splits as $\psi$ is varied to a non-zero value. The half-pitchfork plot for $\psi=0$, separates into two curves for $\psi=\frac{\pi}{12}$, resembling an imperfect system. The parameter $\psi$ acts as an \emph{imperfection} parameter perturbing the the perfect system obtained for $\psi=0$. Interestingly, the equilibria chosen as starting solutions lie on branches of solutions continued from other initial solutions. In these instances, the change in the indices is consistent with the observed folds. If a stable equilibrium exists, it transitions to the other stable equilibrium lying on the other curve, when the parameter is varied past the folds. Based on this observation, we illustrate a few potential configurations that may result from the snapping motion. In each case, multiple snapping configurations are possible, and the exact configuration to which the system reaches depends on the energy released during snapping and the associated dissipation.

Moreover, this response is not reproducible solely by varying $\epsilon$,  i.e., it cannot be replicated just by increasing or decreasing $\epsilon$. When $\epsilon$ is varied alone, the system remains on the lower curve. However, snap instability can be reproduced by tuning the combination of $\epsilon$ and $\psi$. Controlling  $\epsilon$ in the vicinity of folds can be easily implemented in mechanical systems. For example, a switching device can be activated by adjusting $\epsilon$ through a linear actuator or thermal expansion.

\subsubsection{Rotating Load}
\label{sssec:Example3}

\begin{figure}[t]
     \centering
     \begin{subfigure}{0.85\textwidth}
         \centering
         \includegraphics[width=1.0\textwidth]{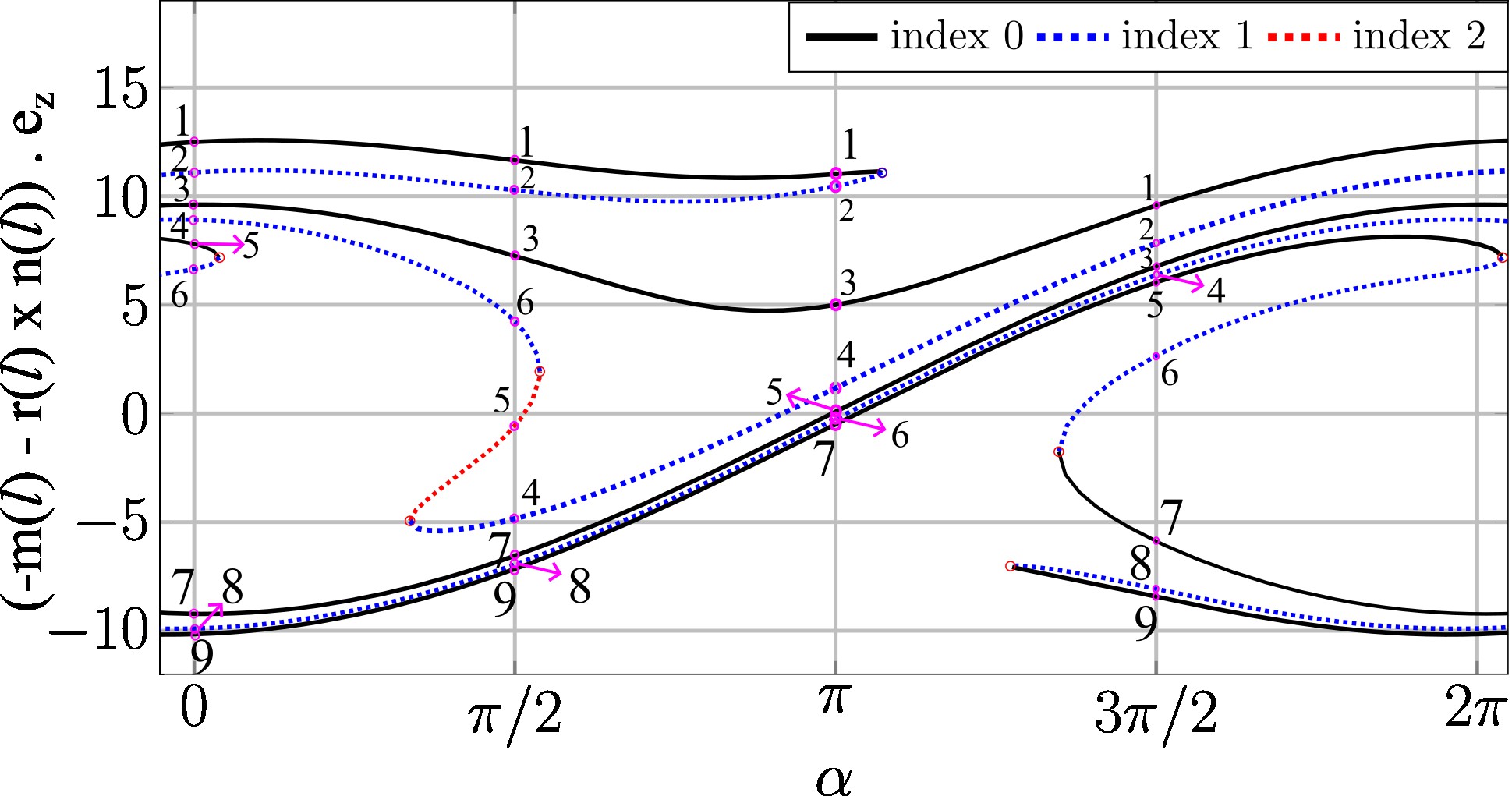}
         \caption{}
         \label{fig:BD}
     \end{subfigure}
     \begin{subfigure}{0.9\textwidth}
         \centering
         \includegraphics[width=1.0\textwidth]{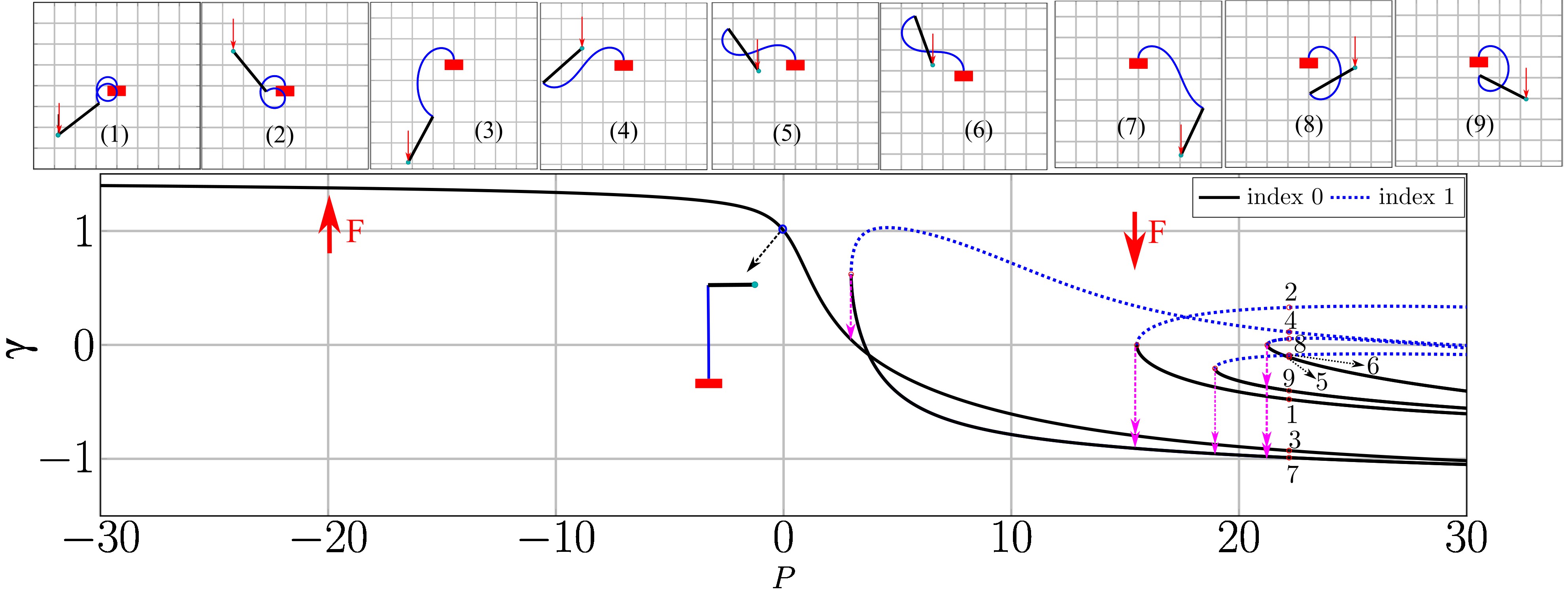}
         \caption{$\alpha=0$}
         \label{fig:Example4_0}
     \end{subfigure}
     \caption{(a) The bifurcation diagram from the previous analysis, adjusted between $0$ and $2 \pi$. The equilibria are labeled based on the order in which the continuation crosses a specified value of $\alpha$.  Continuation is then performed from these labeled equilibria along $P$. (b) The bifurcation diagram when continuation is performed from the equilibria corresponding to $\alpha=0$. The labeled equilibria are displayed both at the top of the plot.}
     \label{fig:Example4_1}
\end{figure}
\begin{figure}[t!]
     \centering
          \begin{subfigure}{0.9\textwidth}
         \centering
         \includegraphics[width=1.0\textwidth]{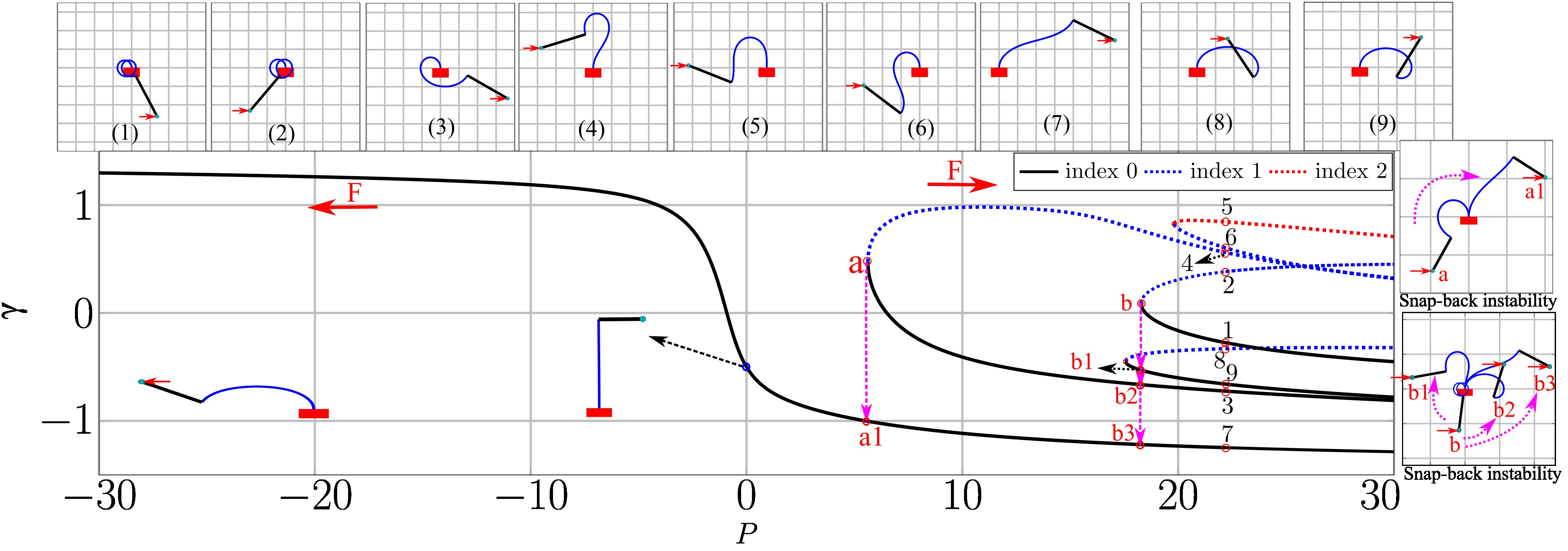}
         \caption{$\alpha=\frac{\pi}{2}$}
         \label{fig:Example4_a}
     \end{subfigure}
    \begin{subfigure}{0.9\textwidth}
         \centering
         \includegraphics[width=1.0\textwidth]{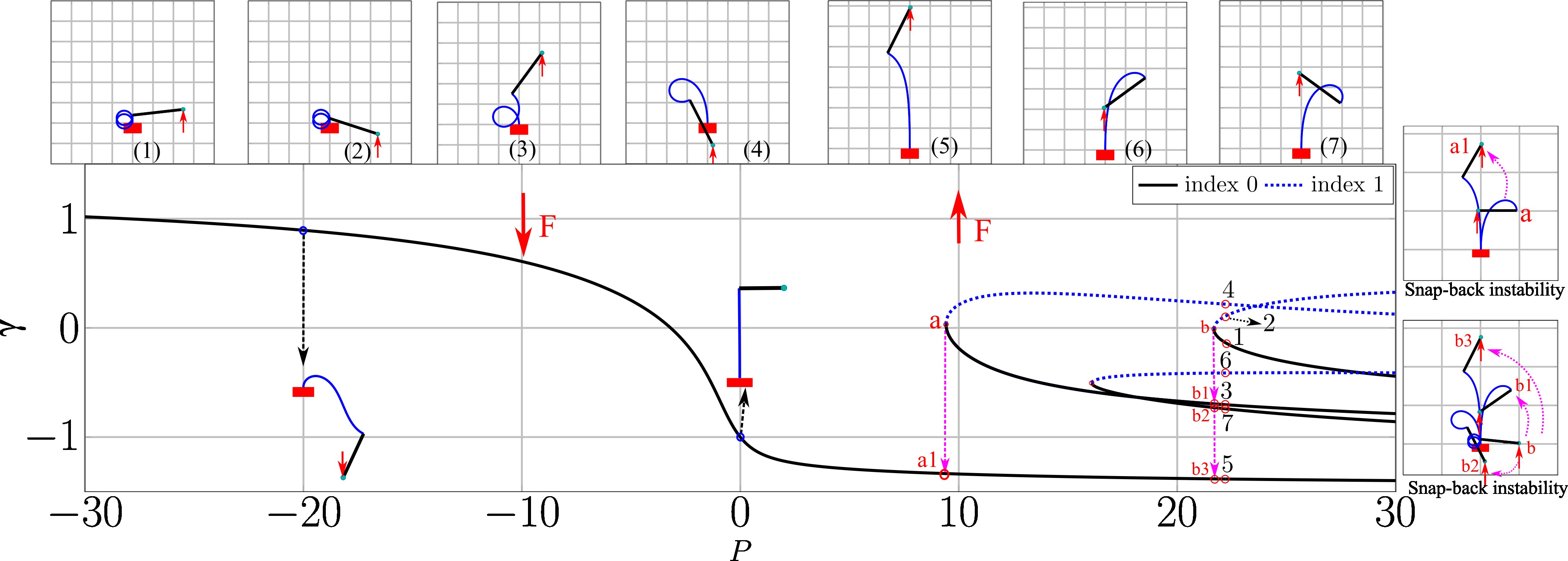}
         \caption{$\alpha=\pi$}
         \label{fig:Example4_b}
     \end{subfigure}
    \begin{subfigure}{0.9\textwidth}
         \centering
         \includegraphics[width=1.0\textwidth]{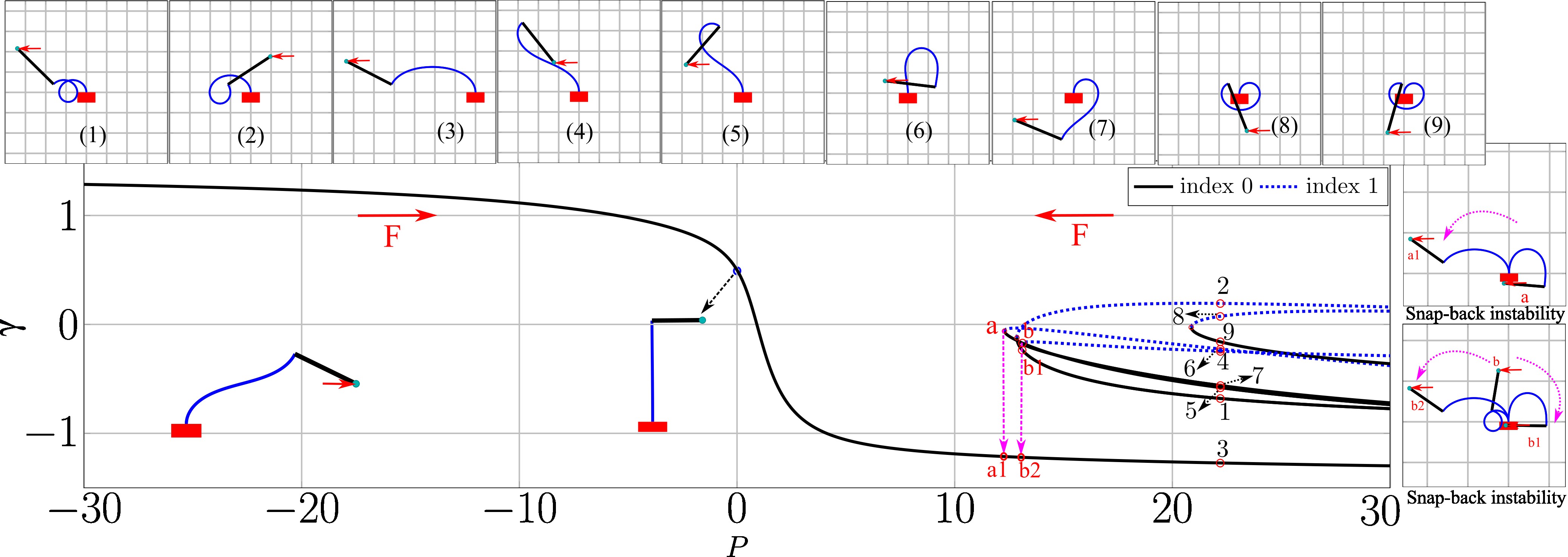}
         \caption{$\alpha=\frac{3\pi}{2}$}
         \label{fig:Example4_c}
     \end{subfigure}
        \caption{ The bifurcation diagram when continuation is performed from the equilibria corresponding to $\alpha=0$. The labeled equilibria and the snap-back instability occurring between the equilibria are shown at the top and on the side.}
        \label{fig:Example4}
\end{figure}

Next, we analyze the dependence of the stability characteristics on the direction of the load $\alpha$. Continuation is performed in $\alpha$ from the equilibrium at $P=\pi^{2}/4$, $\epsilon=0.5$ and $\psi=\pi/2$ as the initial solution. Throughout this subsection, we start from an equilibrium that has an index zero. The ordinate of the bifurcation diagrams in this scenario is 
\begin{align}
\begin{split}
\frac{\partial B}{\partial \alpha} &=-\frac{\partial}{\partial \alpha}\left( \left(P \sin \alpha  \mathbf{e}_{x} - P \cos \alpha \mathbf{e}_{y} \right)\cdot \left( \mathbf{r}(l) + \boldsymbol \Psi (l)\right) \right), \\
&=   \left(\left( \mathbf{r}(l) + \boldsymbol \Psi (l)\right) \times  \mathbf{n}(l) \right) \cdot \mathbf{e}_{z},\\
& =  \left(\mathbf{m}(l) + \mathbf{r}(l) \times \mathbf{n}(l) \right)  \cdot \mathbf{e}_{z}.
\end{split}
\end{align}
and is plotted against $\alpha$ as shown in Figure~\ref{fig:Example3}. A family of equilibria with index zero interspersed with a family of equilibria with index one, is obtained. The elastic configurations before and after the snap-back instability are displayed, illustrating the drastic change. The response of the elastica system reveals additional interesting features. Multiple equilibria exist for any given value of $\alpha$, and configurations are displayed for intermediate values of $\alpha$. In some instances, up to five equilibria can be observed. The self-contact within the present elastica model, as well as contact with the lever arm, are disregarded, and the computed solutions remain valid. This load-direction varying systems in real scenarios can be observed when a rotating magnetic, electric, or gravitational field is applied relative to a fixed frame of elastica. In case of an electric field, a charged particle is attached to the lever arm, while in the case of magnetic field, a ferromagnetic bead is attached to the lever arm.

\subsubsection{Varying Load Magnitude}
\label{sssec:Example4}

We now examine the effect of load's magnitude $P$ on the stability of our elastica setup. We employ the configurations from the previous example as starting solutions and perform continuation along $P$ from $9 \pi^{2}/4$ to $-9 \pi^{2}/4$. The ordinate of the corresponding distinguished bifurcation diagrams is
\begin{align}
\begin{split}
\gamma:= \frac{\partial B}{\partial P}&=-\left( \sin \alpha \mathbf{e}_{x} -\cos \alpha \mathbf{e}_{y} \right) \cdot \left( \mathbf{r}(l) + \epsilon \cos \psi \mathbf{d}_{t}(l)  + \epsilon \sin \psi \mathbf{d}_{r}(l) \right)\\&= -x(l)\sin \alpha + y(l)\cos \alpha + \epsilon \cos \left(\alpha+\psi- \theta(l)\right)  ,
\end{split}
\end{align}
which is the height of the point of application of force projected on the load vector $\mathbf{F}$. Figure~\ref{fig:BD} displays the bifurcation plot from the previous analysis (subsection~\ref{sssec:Example3} Figure~\ref{fig:Example3d}), truncated between $0$ and $2\pi$. Since the system is $2\pi$-periodic with respect to $\alpha$, the section of the plot between $2\pi$ and $4\pi$ is displayed to include all possible equilibria. Note that the section between $0$ and $2\pi$ includes only the equilibria associated with forward continuation. Since backward continuation is not performed in this analysis, the corresponding equilibria are absent. Figure~\ref{fig:Example4_0} and Figure~\ref{fig:Example4} displays the response of the elastica system for this maneuver through bifurcation plots. The equilibria chosen as starting solutions lie on the continuation branches that originate from other initial solutions. The changes in the index are consistent with the predictions (Figure~\ref{fig:schematic_bifurcation_diagram_b}). These diagrams also resemble the unfolding of perfect systems. These equilibrium paths collapse onto the equilibrium path of the perfect system, shown in Figure~\ref{fig:Example2_a}b, as the parameters $\psi$ and $\alpha$ approach zero.  We can also notice the response of the branches emanating from the solutions with index $2$. The response is not reversible and reproducible with solely respect to $P$. However, a combination of parameters $\psi$, $\alpha$ and $P$ can be employed to replicate the behavior. This scenario can be employed in switches, where load fields due to gravity, electric field or magnetic field, can be tuned.

\subsection{Varying Parameter in Fixed end }
Finally, we focus on the stability transitions when the varying parameter appears at the fixed end $s=0$. In this case, we have only one parameter that is the rotation of the clamped end $\theta_{o}$.

\subsubsection{Rotating base}
\label{sssec:Example5}
The analysis is performed for a non-zero arm $\Delta=0.5$ held at $\psi=\pi/2$ for two loads $P= \pi^{2}/4$ and $9 \pi^{2}/4$. In this case, the ordinate of the bifurcation diagrams given by~\eqref{eqn:distinguished_biff_ordinate_fixed} is
\begin{align*}
\left[ \frac{\partial \mathcal{L}}{\partial \zeta^{\prime}} \cdot \frac{\partial \zeta}{\partial \theta_{o}} \right]_{s=0} &= K \theta^{\prime}(0) \\&\equiv \mathbf{m}(0)\cdot \mathbf{e}_{z},
\end{align*}

\begin{figure}[t!]
     \centering
     \begin{subfigure}{\textwidth}
         \centering
        \includegraphics[width=0.75\textwidth]{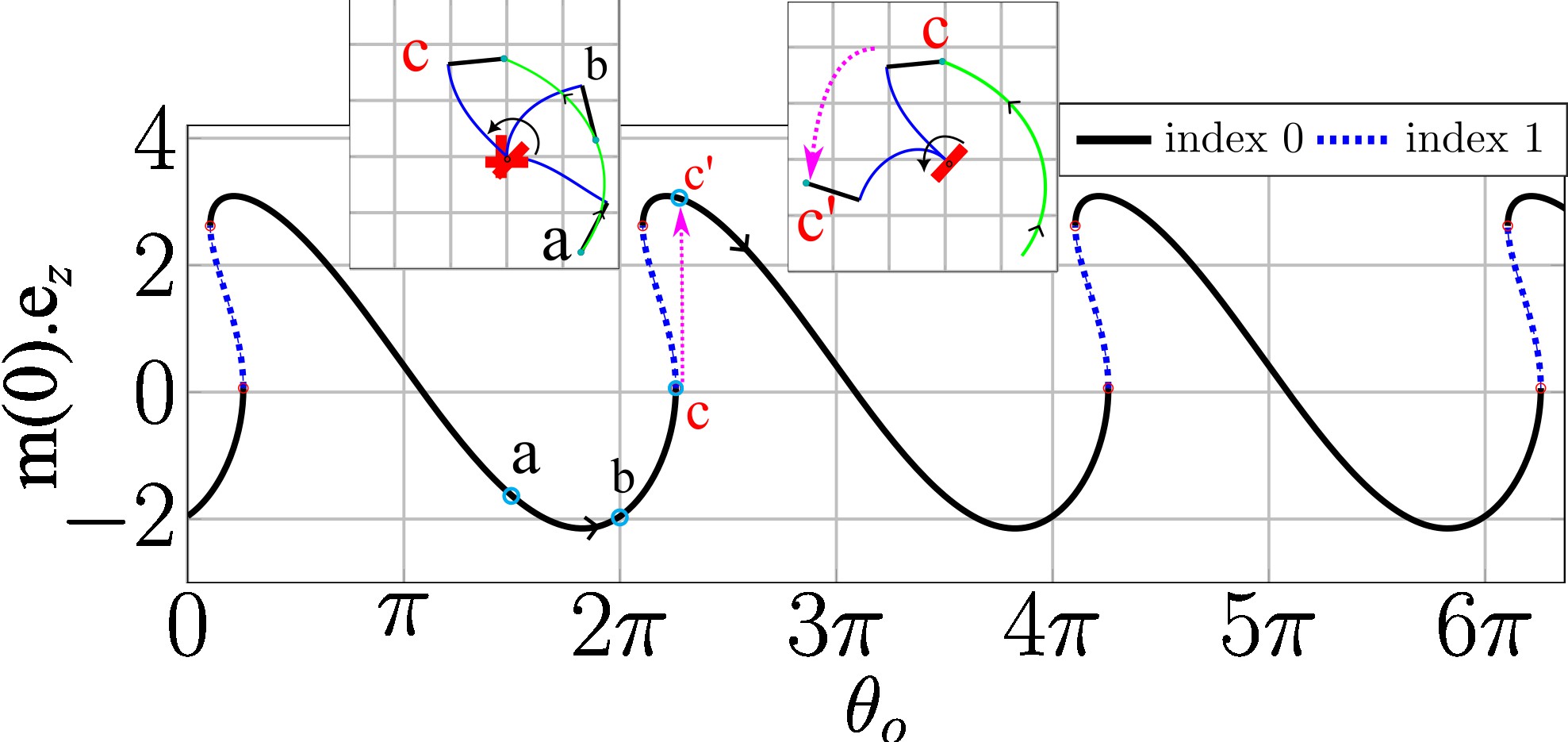}
         \caption{$P = \frac{\pi^{2}}{4}$}
         \label{fig:Example5a}
     \end{subfigure}
     \begin{subfigure}{\textwidth}
         \centering
         \includegraphics[width=0.75\textwidth]{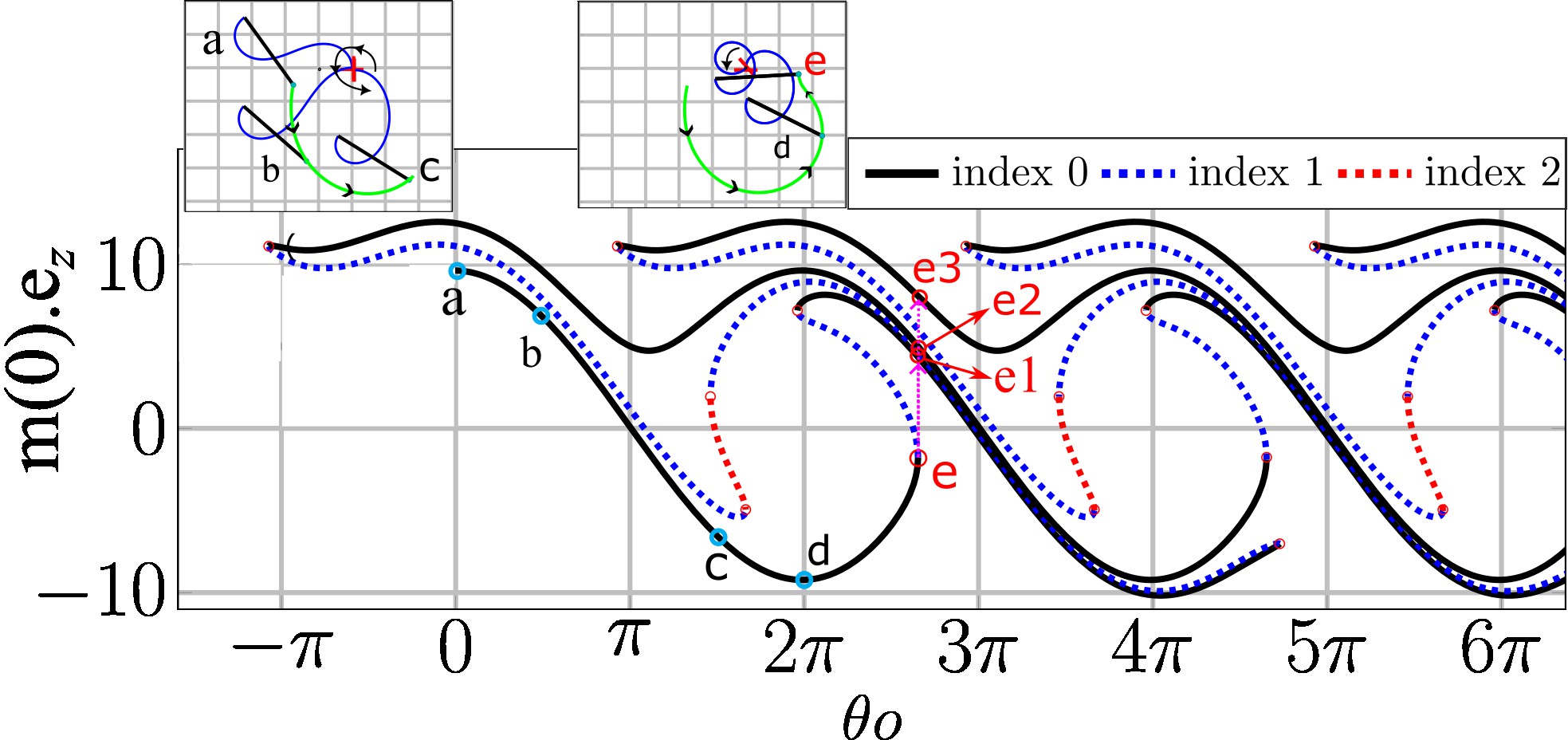}
         \caption{$P =\frac{9 \pi^{2}}{4}$}
         \label{fig:Example5b}
     \end{subfigure}
        \caption{The bifurcation plots as the clamped end $\theta_{o}$ is varied. Selected cases of stable equilibria are labeled and displayed along with the potential snap-back instability. The snap-back instability cases for $P = \frac{9 \pi^{2}}{4}$ are shown in Figure~\ref{fig:Example5_2}. The tip trace before the snap-back instability is shown in green.}
        \label{fig:Example5}
\end{figure}

Figure~\ref{fig:Example5} displays the bifurcation plots for loads $P=\pi^{2}/4$ and $P=9\pi^{2}/4$. The presence of folds indicates an exchange of stability, and the direction of these changes aligns well with the predictions of the bifurcation diagrams (Figure~\ref{fig:schematic_bifurcation_diagram_a}). Indeed, this represents the simplest scenario that can be realized in practical devices. Physically, this problem is equivalent to the system with rotating loads discussed in subsection~\ref{sssec:Example3} of section~\ref{ssec:Examples_free_end}. However, the ordinates in the distinguished bifurcation diagrams differ. Nevertheless,  these ordinates prove to be equivalent when~\eqref{eqn:Newton_balance} are examined, which essentially represents that the quantity
\begin{align*}
\mathbf{m}+ \mathbf{r} \times \mathbf{n} 
\end{align*}
is constant for $s \in [0,l]$. Therefore,
\begin{align*}
\mathbf{m}(0)=\mathbf{m}(l) +\mathbf{r}(l) \times \mathbf{n}(l) .
\end{align*}

\begin{figure}
         \centering
         \includegraphics[width=0.75\textwidth]{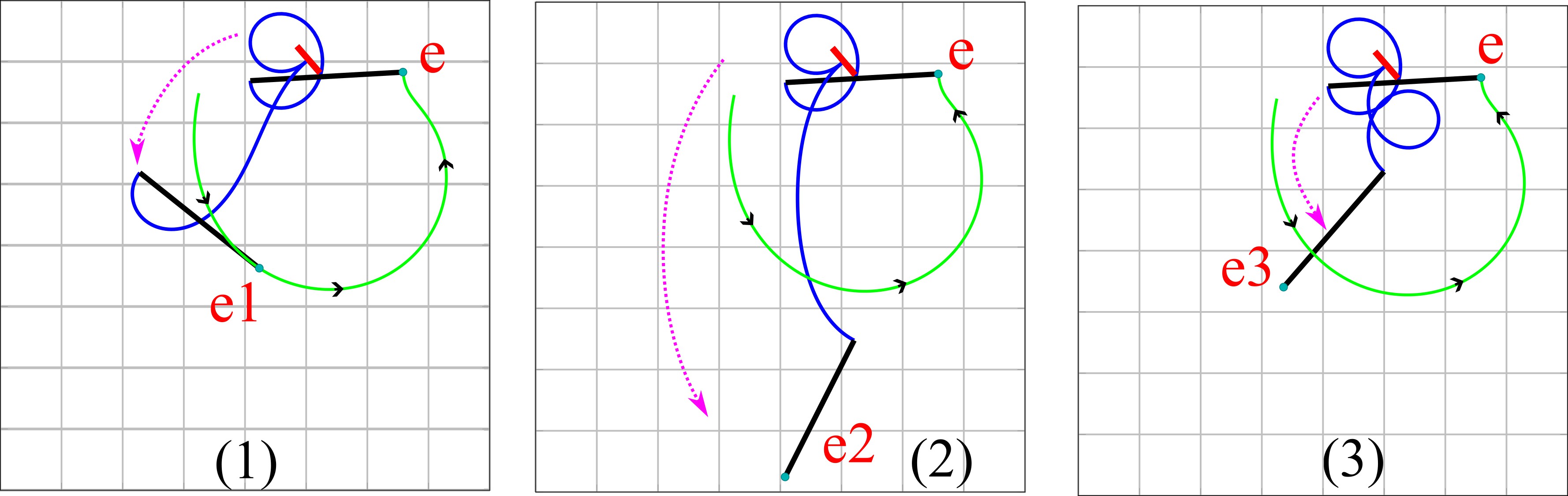}
         \caption{Configurations before and after the snapping for the case shown in Figure~\ref{fig:Example5b}. The tip trace prior to the snap-back instability is displayed in green.}
         \label{fig:Example5_2}
\end{figure}

In our setup, the parameters $\alpha$ and $\theta_{o}$ are interdependent; varying $\alpha$ is equivalent to varying $\theta_{o}$ in the negative direction. As a result, the ordinate in the current case is same as that presented in  subsection~\ref{sssec:Example3} of section~\ref{ssec:Examples_free_end}, while the abscissa is sign-reversed. The plots in Figure~\ref{fig:Example5} represent horizontally mirrored versions of those shown in Figure~\ref{fig:Example3}. The evolution of elastica configurations as the clamped end is rotated, are also indicated in the plots, along with the snap-back instability, illustrating the catapult-like behavior. The self-contact within the elastica and with the lever arm is disregarded. The lower load case of $\pi^{2}/4$ led to equilibria that don't involve self-contact. The higher load case of $9\pi^{2}/4$ generated multi-stable equilibria that formed loops and the snap-back instability would result in more than one equilibrium. Two families of stable equilibrium can be observed. The equilibrium configurations of one family are illustrated in Figure~\ref{fig:Example5}. The configurations for the second family, along with the snap-back instability, are shown in Figure~\ref{fig:Example5_2}. 

\section{Conclusion}
\label{sec:Summary}
The theory of distinguished bifurcation diagrams was extended to problems characterized by fixed-free ends for their stability analysis. Using a combination of bifurcation diagrams and fold information, the stability of elastica subjected to an end load through a rigid lever arm was studied. The stability is determined through a qualitative examination of the plots, without performing a rigorous analysis such as conjugate point computations or eigenvalue  determination. Nevertheless, the stability index of at least one solution along a branch must be established, which was determined using bifurcation analysis. Despite their simplicity and relative abundance, elastica systems with a lever arm have received limited attention. We presented several instances of multi-stability arising in these special elasticas by tuning various parameters. Bi-stability or multi-stability is extensively investigated for engineering applications such as microfluidic devices~\cite{Hilber2016}, soft robots~\cite{Laschi2012, MajidiCarmel2014} and MEMS devices~\cite{zhang2007snap, maurini2007distributed}. Many modern functional materials~\cite{shen2020stimuli} have been developed, where electric, magnetic, optical, thermal, or solvent-based stimuli can be used to tune system parameters near instability, and thereby enabling actuation. The presented elastica system holds significant potential for designing variety of innovative mechanisms. The distinguished bifurcation diagrams can serve as an invaluable tool in their design. A natural extension of this work would be to apply a similar framework to three-dimensional elastic rods using the Kirchhoff rod theory~\cite{Dichmann1996}. 

We restricted our analysis to parameters that produce equilibria with a maximum index of $2$. Future work could explore parameters that generate higher-index equilibria. For example, larger values of $P$ would generate equilibria with higher indices. Another key aspect of the current findings is the system's periodic nature with respect to the rotation parameters. In contrast to problems with fixed-fixed ends in \cite{Hoffman2005}, the present system exhibited periodicity with respect to the rotation parameters when one of its fixed ends is set free. The free end allows for more freedom, limiting the possibility of configurations with higher stored elastic energy. But, increasing the parameters, such as $P$ allows their production.

In this study, we presented the equilibria resulting from the quasi-static control of parameters and drew conclusions on the equilibria after snap-back instability without detailing the transition path. In some instances, multiple stable equilibria are possible after the snap-back instability, and exact equilibrium to which the system transitions cannot be concluded. In this case, the study of dynamical aspects of snapping would be more advantageous~\cite{Armanini2017,snyder1990dynamics}. The most likely equilibrium depends on the energy released during snapping as well as the system's dissipation. Incorporating isoperimetric constraints, where the position of the free end is fixed, while it is free to rotate, would lead to interesting scenarios. It has already been proved that these constraints would have no effect on the ordinate of bifurcation diagrams for the case of fixed-fixed ends~\cite{Hoffman2005}, and it remains to be verified whether the same holds in the present case. The current analysis focuses solely on simple folds. An extension to degenerate scenarios involving non-simple folds would be a valuable direction for future work.

%The current research findings demonstrating the several instances of multi-stability have implications in applications in Engineering and technology such as in the design of novel soft robotic arm mechanisms, switching/triggering mechanisms or miniature MEMS devices. 

\section*{Acknowledgments}
I thank Prof. John Maddocks for fruitful discussions and for sharing his extensive knowledge of variational principles and elastic rods. This work would not have been possible without his guidance. I also thank Prof. Raushan Singh and Prof. Harmeet Singh for their insightful comments on the initial drafts of this work. This study was funded by the Einstein Foundation Berlin and the Deutsche Forschungsgemeinschaft (DFG, German Research Foundation) under Germany's Excellence Strategy – The Berlin Mathematics Research Center MATH+ (EXC-2046/1, project ID: 390685689). I also acknowledge funding from the Italian Ministry of University and Research in the framework of the Call for Proposals for scrolling of final rankings of the PRIN 2022 call - Protocol no. 2022WFJ795.

%Bibliography
\bibliographystyle{unsrt}  
\bibliography{Arxiv}

\end{document}